\documentclass[10pt,journal]{IEEEtran}

\ifCLASSINFOpdf
\else
\fi
\usepackage[cmex10]{amsmath}
\usepackage{amssymb}
\usepackage{enumitem}
\usepackage{subfigure}
\usepackage{graphicx}

\newcommand{\kron}{\raisebox{1pt}{\ensuremath{\:\otimes\:}}} 
\DeclareMathOperator*{\vech}{vec}
\DeclareMathOperator*{\Var}{Var}
\newcommand{\bz}{\mathbf{z}}

\newcommand{\R}{\mathbb{R}}
\newcommand{\C}{\mathbb{C}}
\newcommand{\Z}{\mathbb{Z}}
\newcommand{\N}{\mathbb{N}}

\usepackage{amsthm}
\newtheorem{theorem}{Theorem}
\newtheorem{cor}{Corollary}

\renewcommand{\IEEEQED}{\IEEEQEDopen}
\usepackage{color}

\begin{document}
%
\title{Validity of time reversal for testing Granger~causality}
%
%
%

\author{Irene Winkler, Danny Panknin, Daniel Bartz, Klaus-Robert M\"uller and Stefan Haufe
\thanks{This work was supported by a Marie Curie International Outgoing Fellowship (grant No. 625991) within the 7th European Community Framework Program, the BMBF project ALICE II, Autonomous Learning in Complex Environments (01IB15001B), and the Brain Korea 21 Plus Program as well as the SGER Grant 2014055911 through the National Research Foundation of Korea funded by the Ministry of Education. }
\thanks{I. Winkler, D. Panknin, D. Bartz, K.-R. M\"uller and S. Haufe are with the Machine Learning Group, Technische Universit\"at Berlin, Germany. K.-R. M\"uller is also with the Department of Brain and Cognitive Engineering, Korea University, Seoul, Republic of Korea. S. Haufe is also with the Laboratory for Intelligent Imaging and Neural Computing, Columbia University, New York, USA.
Correspondence to: \{i.winkler,stefan.haufe\}@tu-berlin.de.}
}

%
%

\markboth{Validity of time reversal for testing Granger causality}{}
%



\maketitle

\begin{abstract}
Inferring causal interactions from observed data is a challenging problem, especially in the presence of measurement noise. To alleviate the problem of spurious causality, Haufe et al. (2013) proposed to contrast measures of information flow obtained on the original data against the same measures obtained on time-reversed data. They show that this procedure, time-reversed Granger causality (TRGC), robustly rejects causal interpretations on mixtures of independent signals. While promising results have been achieved in simulations, it was so far unknown whether time reversal leads to valid measures of information flow in the presence of true interaction. 
Here we prove that, for linear finite-order autoregressive processes with unidirectional information flow between two variables, the application of time reversal for testing Granger causality indeed leads to correct estimates of information flow and its directionality. Using simulations, we further show that TRGC is able to infer correct directionality with similar statistical power as the net Granger causality between two variables, while being much more robust to the presence of measurement noise. 
\end{abstract}

\begin{IEEEkeywords}
Granger causality, time reversal, noise, TRGC
\end{IEEEkeywords}

%
\IEEEpeerreviewmaketitle

\section{Introduction}
\label{intro}
%
%
%
%

\IEEEPARstart{T}{he} estimation of causal relations between time series is a signal processing topic promising to enhance our understanding of dynamical systems in numerous application domains. For data with time structure, the concept of Granger causality (GC) has gained popularity as a simple testable definition of causality based on temporal precedence. Signal processing techniques based on Granger-causality have been studied in a variety of fields such as econometrics~\cite{Luedtkepohl2007}, 
neuroscience \cite{Winterhalder2005,Mader2008,Bressler2011,Bolstad2011}, and climate science \cite{Kaufmann1997,Triacca2001}. 

In its original formulation, a time series $x_t$ is said to Granger-cause a time series $y_t$, if the past of $x_t$ helps to predict $y_t$ above what can be predicted by using `all other information in the universe' besides the past of $x_t$ \cite{Granger1969}. 
In the bivariate framework, it is common to consider only the information contained in the past of $x_t$ and $y_t$ (cf. \cite{Hamilton1994}).

A serious problem for the estimation of information flow using Granger causality is that spurious Granger causality can occur due to measurement noise. On one hand, if two sensors measuring the same signal are superimposed with noise, they mutually help predicting each other's future \cite{Nalatore2007,Nolte2008}. This is a problem especially in the study of brain connectivity using non-invasive electrophysiology, where the activity at a given sensor is typically a mixture of contributions from several neuronal sources due to the volume conduction of electric currents in the head 
\cite{Nolte2004,GomezHerrero2008,Schoffelen2009,Haufe2010,Haufe2013}. On the other hand, noise that is correlated across sensors has a similar adverse effect on estimates of directed interaction even if the actual 
signals-of-interest are not mixed into different sensors \cite{Haufe2012,Vinck2015}. Such spurious causality can occur in any measure based on the concept of Granger causality, including multivariate \cite{Kaminski1991,Baccala2001} and non-linear \cite{Marinazzo2008PhysRevLett,GrosseWentrup2009,Vicente2011} variants.

Recently, a number of ways to make causality estimates more robust to the presence of mixed signals and noise have been proposed. These include novel measures of directed information flow \cite{Nolte2004,Nalatore2007,Nolte2008} as well as novel ways of assessing their statistical significance \cite{Breakspear2004,Vicente2011,Haufe2012,Haufe2013,Vinck2015}. Recently, Haufe et. al \cite{Haufe2012,Haufe2013} suggested to contrast causality scores obtained on the original time series to those obtained on time-reversed signals. The intuitive idea behind this approach is that, if temporal order is crucial to tell a driver from a recipient, directed information flow should be reduced (if not reversed) if the temporal order is reversed. In fact, Haufe et al. showed that for correlated, but non-interacting signals, the use of time reversal for testing Granger causality scores (here referred to as time-reversed Granger causality, TRGC) and other metrics based on cross-spectral estimates or linear autoregressive modeling correctly leads to rejection of causal interpretations. This was confirmed for Granger causality in an independent simulation study \cite{Vinck2015} showing that TRGC leads to a much smaller fraction of false positive detections compared to the original Granger causality index, and also compares favorably against the Phase Slope Index (PSI) \cite{Nolte2008}. 

While time-reversed Granger causality thus displays an intriguing noise robustness property, and yields very encouraging results in simulations, its behavior \emph{in the presence of causal interactions} is still poorly understood. In particular, it is currently unclear how Granger causality scores computed on time-reversed signals link to the causal interactions on the original time-series, and therefore whether TRGC correctly indicates the direction of causality. Theoretical guarantees have only been derived for special cases in which either the signal's auto- and cross-covariances are very small in magnitude, or in which both signals have very similar autocorrelations \cite{Vinck2015}. 

The aim of this paper is two-fold. In the theory section, we provide new theoretical insights on time-reversal for testing Granger causality between two variables. After introducing the concepts of linear autoregressive modeling, Granger causality, and time-reversed Granger causality (Section~\ref{sec:var_gc} and \ref{sec:timereversal}), we elaborate on the existing result of Haufe et al. \cite{Haufe2013} showing that, for mixtures of independent signals, causality measures based on  cross-covariances are invariant to the reversal of the temporal order (Section~\ref{sec:NoiseRobustness}). This is the theoretical basis for the noise-robustness property of time-reversal testing of causality scores.
We then investigate the time-reversal of a process fulfilling the assumptions typically made by Granger causality estimators: a finite-order vector autoregressive (VAR) process that is unaffected by measurement noise. We review what is known about the time-reversal of a VAR process (Section~\ref{sec:timereversedAR}), based on what we provide an analytic description of Granger causality scores of
time-reversed signals in terms of their autoregressive coefficients (Section~\ref{sec:TRGCdescription}) and a minimal example (Section~\ref{sec:minimal}). Using these insights, we prove our main result stating that, in the case of unambiguous unidirectional information flow from $x_t$ to $y_t$, time reversal leads to a decrease of the Granger-causal net information flow 
relative to the original time series. The difference of net Granger causality scores obtained on original and time-reversed data thus indicates the correct direction of interaction (Section~\ref{sec:maintheorem}). 

In the second part of the paper (Section~\ref{sec:experiments}), we revisit scenarios known to cause problems for conventional Granger causality. Using simulations, we illustrate when and how the theoretical guarantees of TRGC  
lead to measurable performance increases in practice. We point out the implications of our theoretical and empirical results in Section~\ref{sec:Discussion}, along with a discussion of ambiguities in causal interpretation caused by the presence of correlated residuals in VAR models. 

\section{Theory}
\label{sec:theory}

Vectors are considered to be column vectors (unless otherwise stated), and are generally typed in bold. The symbol $\cdot^\top$ denotes the transpose operator, $I$ the identity matrix, and $[\cdot, \cdot]$ concatenation. The symbol $\kron$ refers to the Kronecker product, and $\vech (\cdot)$ to the vectorization operator, which converts a matrix into a column vector. The symbol $\langle \cdot \rangle$ denotes expectation. 
The cross-covariance matrices of a stationary process $\bz_t$ are denoted by
\begin{equation*}
C_\bz (h) := \Big\langle (\bz_t - \langle \bz \rangle ) (\bz_{t-h} - \langle \bz \rangle )^\top \Big \rangle \qquad \forall h \in \Z \;.
\end{equation*}
We use the notation $\bz_t$ both for an observed time series and its underlying data generating process. 
We denote all quantities related to the time-reversed process $\tilde{\bz}_t:=\bz_{-t}$ with a tilde.

A process $\boldsymbol{\epsilon}_t$ is said to be \emph{white noise} if it is stationary with mean zero, finite covariance and zero autocorrelation; that is, if $C_\epsilon (h) = 0 \; \forall h \in \Z \setminus \{ 0\}$. Note that the covariance matrix $C_\epsilon(0)$ is not necessarily diagonal, and that neither independence nor joint Gaussianity is required. 

\subsection{Granger causality and the linear VAR model}
\label{sec:var_gc}

Consider a stable bivariate vector autoregressive process of lag order $p$ (VAR($p$) process),  $\bz_t  = \begin{bmatrix} x_t \\ y_t \end{bmatrix} \in \R^2$, 
\begin{equation}
\bz_t =  A_1 \bz_{t-1} + A_2 \bz_{t-2}  + \ldots + A_p \bz_{t-p}  + \epsilon_t  \;,\label{eq:VARp}
\end{equation}
where $\epsilon_t \in \R^2$ is a $2$-dimensional white noise process (that is, $\langle \epsilon_t \rangle = 0$,  $\langle \epsilon_t \epsilon_{t-h}^\top \rangle = 0$ for $h  \in \Z \setminus \{0\}$, and $\langle \epsilon_t \bz_{t-h}^\top \rangle = 0$ for $h \in \N \setminus \{0\}$) with residual covariance matrix
\begin{equation}
\Sigma = \langle \epsilon_t \epsilon_t^\top \rangle = \begin{bmatrix} \Sigma_{xx} & \Sigma_{xy} \\ \Sigma_{xy} & \Sigma_{yy} \end{bmatrix} \;. 
\end{equation}
The noise variables $\epsilon_t$ are also called \emph{innovations} {or \emph{residuals}. Stability requires that $ \det(I - A_1 \lambda - \ldots - A_p \lambda^p) \neq 0$ for all $\lambda \in \C$ with $|\lambda|\leq 1$.}

Following \cite{Geweke82}, $x_t$ and $y_t$ possess themselves autoregressive (AR) representations, which we denote by
\begin{align}
x_t &= \sum_{k = 1}^\infty a_k x_{t-k} +  \xi_t^x ,  \quad \Var( \xi_t^x)  =: \Sigma_{x}  \quad \text{and} \label{eq:arres} \\
y_t &= \sum_{k = 1}^\infty b_k y_{t-k}  +  \xi_t^y, \quad \Var( \xi_t^y)  =: \Sigma_{y} \;.\label{eq:arres_y}
\end{align}

The residuals $\xi_t^x$ and $\xi_t^y$ of these two univariate processes are each serially uncorrelated, but may be correlated with each other at various lags.  Importantly, even though the bivariate autoregressive process~\eqref{eq:VARp} is of finite order, the univariate processes~\eqref{eq:arres} and \eqref{eq:arres_y} are in general of infinite order. We refer to \eqref{eq:VARp} as the \textit{unrestricted} or \textit{full} model, while \eqref{eq:arres} and \eqref{eq:arres_y} contain the \textit{restricted} models. 

Directed Granger-causal information flow is defined based on the so-called \emph{Granger-scores} \cite{Geweke82}
\begin{equation}
\begin{aligned} 
F_{y \rightarrow x} :=  \log \left( \frac{\Sigma_{x}}{\Sigma_{xx}} \right)  \quad  \text{and} \quad F_{x \rightarrow y} :=  \log \left( \frac{\Sigma_y} {\Sigma_{yy}} \right) \;.
\end{aligned}
\label{eq:Geweke_def}
\end{equation}

Granger causality from $x_t$ to $y_t$ implies that information from the past of $x_t$ improve the prediction of the present of $y_t$ compared to what can be predicted by the past of $y_t$ alone. That is, the residual variance $\Sigma_{yy}$ of the unrestricted model is required to be smaller than the residual variance $\Sigma_y$ of the restricted model. Under the assumption of Gaussian-distributed residuals, $F_{y \rightarrow x}$ and $F_{x \rightarrow y}$ are asymptotically $\chi^2$ distributed, giving rise to an analytical test of their significance \cite{Geweke82}. An asymptotically equivalent test is given by an F-test of the goodness-of-fit of the two models (cf. \cite{Hamilton1994,Bressler2011}). We refer to this approach as \emph{standard Granger causality} (standard GC).

As variables in physical systems often mutually influence each other, it is also of interest to determine the \emph{net} driver of the interaction by assessing whether more information is flowing from $x_t$ to $y_t$ then from $y_t$ to $x_t$ or vice versa. Following \cite{Nolte2008,Nolte2010}, \emph{net Granger causality} (Net-GC) is defined as the difference of the Granger causality scores, that is  
\begin{align}
F_{x \rightarrow y}^{(net)}  := F_{x \rightarrow y} - F_{y \rightarrow x} \quad \text{and} \quad F_{y \rightarrow x}^{(net)}  := -F_{x \rightarrow y}^{(net)} \;. \label{eq:Net-GCdef}
\end{align}

As the analytical distributions of these differences are unknown, statistical significance of Net-GC scores needs to be assessed using resampling methods as outlined in Section~\ref{sec:exp_methods}.

\subsection{Time-reversed Granger causality (TRGC)}
\label{sec:timereversal}

To avoid false detections of causal interactions, Haufe et al. proposed to contrast causality measures applied to the original time series with the same measures obtained from time-reversed signals $\tilde{\bz}_t :=\bz_{-t}$ \cite{Haufe2012,Haufe2013}. Here, we formalize this idea in the context of Granger causality.

Given a bivariate VAR($p$) process, its time-reversed process $\tilde{\bz}_t $ also possesses a VAR($p$) representation, which we derive in Section~\ref{sec:timereversedAR}. We denote the residual covariance matrix of the time-reversed process by
\begin{equation}
\tilde{\Sigma} = \begin{bmatrix} \tilde{\Sigma}_{xx} & \tilde{\Sigma}_{xy} \\ \tilde{\Sigma}_{xy} & \tilde{\Sigma}_{yy} \end{bmatrix} \;.
\end{equation}

The restricted AR models of the time-reversed data have a simple structure, as they are concerned with univariate time series. 
{The autocovariance function of a univariate time series is symmetric, i.e., we have \mbox{$C_{x}(h) = C_{x}(-h)$} and \mbox{$C_{y}(h) = C_{y}(-h)$} for all $h \in \Z$. As a result of this and~\eqref{eq:reversedCov} (Section~\ref{sec:NoiseRobustness}), the time-reversed signals will have the same autocovariances as the original series.}
Because the AR representation is uniquely determined by the autocovariance function (cf. Section~\ref{subsec:VARCovs}), they also share the same AR representation.
The restricted models of the time-reversed univariate processes are thus given by
\begin{align}
x_t &= \sum_{k = 1}^\infty a_k x_{t+k} +  \tilde{\xi}_t^x ,  \quad \Var( \tilde{\xi}_t^x)  =: \tilde{\Sigma}_{x} \quad \text{and} \quad \label{eq:arresrev} \\
y_t &= \sum_{k = 1}^\infty b_k y_{t+k}  +  \tilde{\xi}_t^y, \quad  \Var(\tilde{\xi}_t^y)  =: \tilde{\Sigma}_{y} \label{eq:aryresrev}
\end{align}
with
\begin{equation}
\tilde{\Sigma}_{x}  = \Sigma_{x} \qquad \mbox{ and } \qquad \tilde{\Sigma}_{y}  = \Sigma_{y} \;. \label{eq:uniequi}
\end{equation}

In analogy to the original time series, we define the time-reversed Granger scores as 
\begin{equation}
\tilde{F}_{\tilde{y} \rightarrow \tilde{x}} :=  \log \left( \frac{\tilde{\Sigma}_{x}}{\tilde{\Sigma}_{xx}} \right)  \quad \text{and} \quad \tilde{F}_{\tilde{x} \rightarrow \tilde{y}} :=  \log \left( \frac{\tilde{\Sigma}_y}{\tilde{\Sigma}_{yy}} \right) \;,
\end{equation}
and the net Granger causality scores as
\begin{align}
\tilde{F}_{\tilde{x} \rightarrow \tilde{y}}^{(net)}  &:=\tilde{F}_{\tilde{x} \rightarrow \tilde{y}} - \tilde{F}_{\tilde{y} \rightarrow \tilde{x}}  \quad \text{and} \quad \tilde{F}_{\tilde{y} \rightarrow \tilde{x}}^{(net)}  := -\tilde{F}_{\tilde{x} \rightarrow \tilde{y}}^{(net)} \;.
\end{align}

Finally, the differences of the Granger scores obtained on original and time-reversed signals are given by
\begin{align}
\tilde{D}_{y \rightarrow x}  & := F_{y \rightarrow x} - \tilde{F}_{\tilde{y} \rightarrow \tilde{x}} \;, \label{eq:D_def1} \\
\tilde{D}_{x \rightarrow y}  & :=  F_{x \rightarrow y} - \tilde{F}_{\tilde{x} \rightarrow \tilde{y}}  \;, \quad \text{and}    \\
\tilde{D}_{x \rightarrow y}^{(net)} & :=  F_{x \rightarrow y}^{(net)}  - \tilde{F}_{\tilde{x} \rightarrow \tilde{y}}^{(net)} \;.
\label{eq:D_def}
\end{align}

Time-reversed Granger causality can be applied in the following variants.

\paragraph{Conjunction-based time-reversed Granger causality (Conj-TRGC)}  
Here, net information flow from $x_t$ to $y_t$ is inferred if 
\begin{equation}
F_{x \rightarrow y}^{(net)} > 0 \quad \textnormal{and} \quad \tilde{F}_{\tilde{x} \rightarrow \tilde{y}}^{(net)} < 0  \;, \label{eq:conjTRGC}
\end{equation}
that is, if the directionality of net Granger causality reverses for time-reversed signals. This variant has been investigated in \cite{Vinck2015}. 

\paragraph{Difference-based time-reversed Granger causality (Diff-TRGC)}  
Here, net information flow from $x_t$ to $y_t$ is inferred if
\begin{equation}
\tilde{D}_{x \rightarrow y}^{(net)} > 0 \;, \label{eq:diffTRGC}
\end{equation}
that is, we require that net Granger causality from $x_t$ to $y_t$ is reduced on the time-reversed signals. Note that this is a weaker requirement than conjunction-based TRGC, as all signals for which \eqref{eq:conjTRGC} holds also fulfill \eqref{eq:diffTRGC}.

\paragraph{Conjunction of Net-GC and Diff-TRGC}  
Finally, we can require both the time-reversed net difference and the net Granger score to be significantly larger than zero in order to infer net information flow from $x_t$ to $y_t$, that is
\begin{equation}
\tilde{D}_{x \rightarrow y}^{(net)} > 0 \quad \textnormal{and} \quad F_{x \rightarrow y}^{(net)} < 0  \;. \label{eq:netdiffTRGC}
\end{equation}

Just as for Net-GC, statistical significance of Conj-TRGC and Diff-TRGC, as well as the combination of Net-GC and Diff-TRGC can be assessed using resampling techniques (see Section~\ref{sec:exp_methods}). 

\subsection{Robustness of time-reversed Granger causality (TRGC)}
\label{sec:NoiseRobustness}

In \cite{Haufe2013} it is pointed out that time-reversed Granger causality robustly rejects causal interpretations for mixtures of non-interacting signals such as correlated noise sources. The mathematical basis for this noise robustness property is the fact that the cross-covariance matrices $\tilde{C}_{\tilde{\bz}}(\cdot)$ of the time-reversed signals are equal to the transposed cross-covariance matrices of the original signals, that is
\begin{equation}
\tilde{C}_{\tilde{\bz}}(h) = \langle \tilde{\bz}_t \tilde{\bz}_{t-h}^\top \rangle = \langle \bz_t \bz_{t+h}^\top \rangle = C_{\bz}(-h) = \big(C_{\bz}(h)\big)^\top \label{eq:reversedCov}
\end{equation}
for all $h \in \Z$. If a series $\eta_t$ only contains a mixture of independent signals, all its cross-covariance matrices are symmetric \cite{Nolte2006}: consider $\eta_t = M s_t$ where $s_t$ contains a number of independent sources. Then, for all $h \in \Z$, $C_s(h) = diag$ and
thus $C_{\eta}(h) = M C_{s} (h) M^\top$ 
is symmetric. 
For mixtures of independent noise sources, any causality measure that is solely based on a series' cross-covariance matrices therefore yields the same result on the original and the time-reversed signals. This includes Granger causality, but also other popular variants such as directed transfer function (DTF) \cite{Kaminski1991} and partial directed coherence (PDC) \cite{Baccala2001}. Given sufficient amounts of data, the conditions for Conj-TRGC and Diff-TRGC cannot be  fulfilled for mixtures of independent sources using these measures, preventing the detection of spurious interaction.

\subsection{The VAR representation of a time-reversed process}
\label{sec:timereversedAR}

There is so far no theoretical argument guaranteeing that time-reversed Granger causality correctly indicates the presence of information flow as well as its direction
\emph{in the presence of actual interaction}. In order to provide such a guarantee, we here study the time-reversal of (linear) finite-order VAR processes. Note that studying this case is sufficient since, as a results of Wold's decomposition theorem, every stationary, purely nondeterministic, process can be approximated well by a finite order VAR process \cite{Wold1938,Luedtkepohl2007}.

We start by briefly revisiting the link between cross-covariance matrices and VAR representation, which we use throughout the paper, in Section~\ref{subsec:VARCovs}. In Sections~\ref{subsec:TRVAR1} and \ref{subsec:TRVARp}, we then review the theoretical result of Andel \cite{Andel1972} stating that the time-reversed signal of any VAR($p$) process has again a VAR($p$) representation that can be expressed analytically in terms of the original process.  As the description for $p>1$ is mathematically involved, we only treat the case $p=1$ in the main paper, while the proof for arbitrary~$p$ is presented in Appendix~\ref{app:proofs}. 

We use these results to provide an analytic description of difference-based TRGC scores in terms of their autoregressive coefficients (Section~\ref{sec:TRGCdescription}), give a minimal example (Section~\ref{sec:minimal}), and prove our main result stating that, in the case of unambiguous unidirectional information flow, difference-based time-reversed Granger causality indeed yields the correct result (Section~\ref{sec:maintheorem}).

\subsubsection{The cross-covariance function of a VAR process}
\label{subsec:VARCovs}

Most of the insights in this paper are based on the direct link between autoregressive coefficient matrices $A_1, \ldots, A_p$ and residual covariance matrices $\Sigma$ on one hand, and cross-covariance matrix $C_\bz(\cdot)$ on the other hand. This link is established by the Yule-Walker equations as follows (see e.g. \cite{Luedtkepohl2007}). For a VAR(1) process 
\begin{equation}
\bz_t = A_1 \bz_{t-1} + \epsilon_t \label{eq:VAR1} \;,
\end{equation}
the Yule-Walker equations read
\begin{align}
C_\bz(0) & =  A_1 \cdot C_\bz(0) \cdot A_1^\top + \Sigma \quad \text{and} \label{eq:YuleWalker1} \\
C_\bz(h) &= A_1 \cdot C_\bz(h-1) \quad (\forall h  \in \N \setminus \{0 \}) \;. \label{eq:YuleWalker_h}
\end{align}
Given $A_1$ and $\Sigma$, the cross-covariances are uniquely determined from~\eqref{eq:YuleWalker1} through
\begin{equation} 
\vech(C_\bz(0)) = (I - A_1 \kron A_1)^{-1} \vech{\Sigma} \;,
\end{equation}  
while higher-order cross-covariances $C_\bz(h)$ can be recursively computed using \eqref{eq:YuleWalker_h}. Conversely, $A_1$ and $\Sigma$ are uniquely determined by the cross-covariances through 
\begin{align}
A_1 &= C_\bz(1) C_\bz(0)^{-1}  \quad \text{and}  \\
\Sigma &= C_\bz(0) - A_1 C_\bz(0) A_1^\top \;.\label{eq:SigmaFromYuleWalker1}
\end{align}

Results on VAR($1$) processes can typically be extended to higher-order VAR($p$) processes by reducing VAR($p$) processes to their VAR($1$) form. The VAR($1$) representation of a VAR($p$) process as well as the Yule-Walker equations for general VAR($p$) processes are provided in Appendix~\ref{app:YuleWalker}.

\subsubsection{The VAR representation of a time-reversed VAR($1$) process}
\label{subsec:TRVAR1}

The time-reversed autoregressive representation of a VAR($1$) process $\bz_t$ has been derived by Bartlett in 1955 \cite{Bartlett1955}.  Suppose we generate an infinite sequence of $\bz_t$  according to the VAR(1) process~\eqref{eq:VAR1}.  
The VAR representation of the \emph{time-reversed} or \emph{backward} process is given by
\begin{equation}
\bz_t = \tilde{A}_1 \bz_{t+1} + \tilde{\epsilon}_{t} \;,
\end{equation}
where 
\begin{equation}
\tilde{A}_1 =  C_\bz (0) \cdot A_1^\top \cdot C_\bz (0)^{-1} \;,\label{eq:revA1A}
\end{equation}
and where the reversed residuals 
$\tilde{\epsilon}_t$ are calculated from $\bz_t$ as 
\begin{equation}
\tilde{\epsilon}_t := \bz_t -  \tilde{A}_1 \bz_{t+1} 
\end{equation}
with residual covariance matrix
\begin{equation}
\tilde{\Sigma} = \langle \tilde{\epsilon}_t \tilde{\epsilon}_t^\top \rangle = C_\bz (0) - C_\bz (0) \cdot A_1^\top \cdot C_\bz (0)^{-1} \cdot A_1 \cdot C_\bz (0) \;. \label{eq:revA1E}
\end{equation}

It is easy to show that the sequence $\tilde{\epsilon}_t$ is indeed white noise, that is  for all $h \in \Z \setminus \{ 0\}$: $\langle \tilde{\epsilon}_t \cdot \tilde{\epsilon}_{t-h}^\top \rangle = 0$ and for all $h \in \N \setminus \{ 0\}$: $\langle \tilde{\epsilon}_t \cdot \bz_{t+h}^\top \rangle = 0$.  

From~\eqref{eq:revA1A}, we see that the time-reversed coefficient matrix $\tilde{A}_1$ is similar to $A_1$, and thus shares some of its properties, notably its eigenvalues, determinant, trace and rank. However, in the context of Granger causality, it is important to note that many properties of $A_1$ do not transfer to $\tilde{A}_1$. In particular, if $A_1$ is triangular, diagonal, or symmetric, this is not generally the case for $\tilde{A}_1$. 

\subsubsection{The VAR representation of a time-reversed VAR($p$) process}
\label{subsec:TRVARp}

The result of Bartlett on the time-reversed VAR($1$) process has been generalized to VAR($p$) processes by Andel in 1972 \cite{Andel1972}, in a paper that received, so far, little attention. Andel showed that any stable VAR($p$) process~\eqref{eq:VARp} has a time-reversed representation 
\begin{equation}
\bz_t = \tilde{A}_1 \bz_{t+1} + \tilde{A}_2 \bz_{t+2} + \ldots + \tilde{A}_p \bz_{t+p} + \tilde{\epsilon}_t  \label{eq:VARprev}
\end{equation}
that is again of order $p$ with uniquely defined autoregressive coefficients $\tilde{A}_1, \ldots \tilde{A}_p$ and residual covariance matrix $\tilde{\Sigma}$. We reproduce this result in Appendix~\ref{app:timereversedVARp}. 
Note that, while we only treat bivariate VAR processes in this paper, the analytic description of the time-reversed VAR process 
holds for processes of arbitrary dimensionality.

\subsection{Analytic description of Diff-TRGC}
\label{sec:TRGCdescription}

Contrasting Granger scores obtained on original with those obtained on time-reversed signals is simplified by the fact that the AR representation of a univariate time series does not depend on the direction of time.  It follows immediately from \eqref{eq:uniequi}, that the differences of the Granger scores related to original and time-reversed data do not depend on the restricted models:
{
\begin{eqnarray}
\begin{aligned}
\tilde{D}_{y \rightarrow x}  &= F_{y \rightarrow x} - \tilde{F}_{\tilde{y} \rightarrow \tilde{x}}   =   \log \tilde{\Sigma}_{xx} - \log  \Sigma_{xx}  \\
\tilde{D}_{x \rightarrow y}  &= F_{x \rightarrow y} - \tilde{F}_{\tilde{x} \rightarrow \tilde{y}}   =  \log \tilde{\Sigma}_{yy} - \log \Sigma_{yy}  \\
\tilde{D}_{x \rightarrow y}^{(net)} & = F_{x \rightarrow y}^{(net)}  - \tilde{F}_{\tilde{x} \rightarrow \tilde{y}}^{(net)}  \\
&= (F_{x \rightarrow y} - F_{y \rightarrow x} ) - (  \tilde{F}_{\tilde{x} \rightarrow \tilde{y}} -  \tilde{F}_{\tilde{y} \rightarrow \tilde{x}})   \\
&= \log \tilde{\Sigma}_{yy} -   \log \tilde{\Sigma}_{xx}  - \log  \Sigma_{yy}  + \log  \Sigma_{xx}  \;.
\end{aligned}
\label{eq:D_fullmodel}
\end{eqnarray} }

The Granger score differences $\tilde{D}_{y \rightarrow x}, \tilde{D}_{x \rightarrow y}$, and $\tilde{D}_{x \rightarrow y}^{(net)}$ thus only depend on the residual covariance matrices of the full models of the original and time-reversed data. 
For the VAR(1) process, these are given in \eqref{eq:SigmaFromYuleWalker1} and \eqref{eq:revA1E}. For VAR($p$) processes, the residual covariance matrices can be obtained through \eqref{eq:SigmafromYuleWalker} and \eqref{eq:AndelSigma} as described in Appendix~\ref{app:YuleWalker} and  \ref{app:timereversedVARp}.

Please note that while \eqref{eq:D_fullmodel} implies that the unrestricted models can be neglected when computing Granger scores differences, we might gain from including them in finite sample settings. We investigate this issue through simulations in Section~\ref{sec:experiments}.

\subsection{A minimal example} 
\label{sec:minimal}

It is not intuitive to see how the residual variance of the time-reversed process, and thus Granger causality, depends on the autoregressive coefficients of the model. 
Interpretation is made difficult by the occurrence of $C_\bz(0)^{-1}$ in \eqref{eq:revA1E}. 

Let us therefore consider the following minimal case: a VAR(1) process $\bz_t$ with $C_\bz(0) = I$. In that case, $C_\bz(h) = A_1^h$ and $C_\bz^\top(h) = (A_1^\top)^h$ for all $h  \in \Z \setminus \{0\}$ (from \eqref{eq:YuleWalker_h}). All asymmetries in the cross-covariance matrices $C_\bz(h)$ are thus due to asymmetries in $A_1$. 

Furthermore, time-reversing the signal leads to transposition of the autoregressive coefficient matrix  $\tilde{A}_1 = A_1^\top$ as a result of \eqref{eq:revA1A}. The residual covariance matrices \eqref{eq:SigmaFromYuleWalker1} and \eqref{eq:revA1E} are now given by 
\begin{equation*} 
\Sigma =  I - A_1 \cdot A_1^\top \quad \text{and} \quad
\tilde{\Sigma} = I - A_1^\top \cdot A_1  \;.
\end{equation*}

Denote with $A_1 = \begin{bmatrix} a_{11} & a_{12} \\ a_{21} & a_{22}  \end{bmatrix}$ the autoregressive coefficients.  We then have 
\begin{equation*} 
\Sigma_{xx} =  1 - a_{11}^2 - a_{12}^2 \quad , \quad
\tilde{\Sigma}_{xx} = 1 - a_{11}^2 - a_{21}^2  \;,
\end{equation*}
and 
\begin{align*} 
\tilde{D}_{y \rightarrow x} 
=  \log \tilde{\Sigma}_{xx}  - \log \Sigma_{xx}  > 0 
& \Leftrightarrow  \Sigma_{xx} < \tilde{\Sigma}_{xx}  \\
& \Leftrightarrow   a_{12}^2  > a_{21}^2  \;.
\end{align*}

The difference of the Granger scores computed on the original and time-reversed time series thus indicates the correct \emph{net} direction of information flow. We will in general not be able to infer whether~$x_t$ has a Granger-causal influence on~$y_t$. However, we will be able to tell whether~$x_t$ Granger-causes~$y_t$ more than~$y_t$ Granger-causes~$x_t$, or vice versa. 

While this simple case will almost never occur in practice, we give theoretical guarantees for more general cases in the next section. 

\subsection{Validity of TRGC for unidirectional information flow}
\label{sec:maintheorem}

We now prove our main result, the validity of difference-based time-reversed Granger causality in the presence of unidirectional information flow. Consider a bivariate VAR($p$) process with unambiguous unidirectional information flow. This is the case when all coefficient matrices are triangular and the residual covariance matrix $\Sigma$ is diagonal. Then the following theorem holds.

\begin{theorem} 
\label{prop:toprove}
Let $\bz_t = \begin{bmatrix} x_t \\ y_t \end{bmatrix} \in \R^2$ be a stable bivariate VAR($p$) process \eqref{eq:VARp} with the time-reversed representation \eqref{eq:VARprev}. 
Under the assumptions
\begin{enumerate}[label={\upshape(A\arabic*)},leftmargin=1cm]
\item $A_1, \ldots, A_p$ are lower triangular matrices (i.e., $x_t$ may Granger-cause $y_t$, but $y_t$ does not Granger-cause $x_t$), and 
\item $\Sigma$ is a diagonal matrix, i.e. $\Sigma_{xy} = 0$ (the residuals are uncorrelated), and 
\item $C_{\bz}(0)$ is invertible \;,
\end{enumerate}
it holds that 
\begin{align} 
\tilde{\Sigma}_{xx} & \leq \Sigma_{xx}   \label{eq:abschaetzung1} \;,
\intertext{and that} 
\tilde{\Sigma}_{yy} & \geq \Sigma_{yy}  \;. \label{eq:abschaetzung2}
\end{align}
\end{theorem}

\begin{cor} 
Under assumptions {\upshape(A1)-(A3)}, Theorem~1 and~\eqref{eq:D_fullmodel} immediately imply  the following inequalities for the differences of Granger scores:
\begin{align}
\tilde{D}_{y \rightarrow x} = F_{y \rightarrow x} - \tilde{F}_{\tilde{y} \rightarrow \tilde{x}}   &\leq 0 \label{eq:Fyxdiff}\\
\tilde{D}_{x \rightarrow y} = F_{x \rightarrow y} - \tilde{F}_{\tilde{x} \rightarrow \tilde{y}}   &\geq 0   \label{eq:Fxydiff} \\
\tilde{D}_{x \rightarrow y}^{(net)} = (F_{x \rightarrow y} - F_{y \rightarrow x} ) - (  \tilde{F}_{\tilde{x} \rightarrow \tilde{y}} -  \tilde{F}_{\tilde{y} \rightarrow \tilde{x}})   &\geq 0  \label{eq:Fnetdiff} \;.
\end{align} 
\end{cor}

As a result of Corollary~1, net Granger-causal information flow from $x_t$ to $y_t$ is reduced or remains the same when the signal is time-reversed. Thus, in the case of unambiguous unidirectional information flow, difference-based time-reversed Granger causality yields the correct result. Note that it is not true in general that the net flow between the time-reversed signals $\tilde{x}_t$ and $\tilde{y}_t$, $\tilde{F}_{\tilde{x} \rightarrow \tilde{y}}^{(net)}$, is negative (reverses compared to the original series).
That is, conjunction-based TRGC might in some cases incorrectly reject the presence of true causal interaction. 

{Corollary 1 states that each of the three difference scores, $\tilde{D}_{y \rightarrow x}$, $\tilde{D}_{x \rightarrow y}$, and $\tilde{D}_{x \rightarrow y}^{(net)} $ alone is sufficient to infer the correct directionality under assumptions (A1)--(A3). As (A1) requires information flow to be unidirectional, the individual scores $\tilde{D}_{y \rightarrow x}$ and $\tilde{D}_{x \rightarrow y}$ only indicate \emph{net} information flow, which is what is also observed in Section~\ref{sec:minimal}.} 

{The three scores will behave differently if the assumption of uncorrelated residuals (A2) is violated. Then, $\tilde{\Sigma}_{xx} \leq \Sigma_{xx}$ and $\tilde{D}_{y \rightarrow x} \leq 0$ still hold, but the inequalities $\tilde{\Sigma}_{yy} \geq \Sigma_{yy}$, $\tilde{D}_{x \rightarrow y} \geq 0$ and $\tilde{D}_{x \rightarrow y}^{(net)} \geq 0$ do not. On average, the net difference $\tilde{D}_{x \rightarrow y}^{(net)}$ (which equals $\tilde{D}_{x \rightarrow y} - \tilde{D}_{y \rightarrow x}$) is less affected by the presence of correlations in the residuals than any of the individual scores, which is why we defined difference-based TRGC based on $\tilde{D}_{x \rightarrow y}^{(net)}$ in \eqref{eq:diffTRGC}. Nevertheless, all three scores are valid measures for net information flow, as residuals should be uncorrelated if the VAR model accurately describes a physical process.}
The significance of uncorrelated as opposed to correlated residuals is discussed in Section~\ref{sec:corrRes}.

\emph{Sketch of the proof.} The first inequality \eqref{eq:abschaetzung1} is relatively easy to prove. The intuition is the following: Since $y_t$ does not Granger-cause $x_t$, the prediction of $x_t$ is only based on past $x_t$. In contrast, the coefficient matrices $\tilde{A}_1$, \ldots, $\tilde{A}_p$ of the time-reversed representation are  in general not triangular. This means that prediction of the time-reversed signals $\tilde{x}_t$ is not only based on past $\tilde{x}_t$, but can also use information from past $\tilde{y}_t$. We would thus expect that $\tilde{x}_t$ can be better predicted than $x_t$, and that the corresponding residuals are smaller.  

The proof of the second inequality \eqref{eq:abschaetzung2} is more involved. The intuition is the following: we would expect that the `amount' of unexplainable variance is the same for both the original and the time-reversed process. Thus, since the residual variance of $x_t$ decreases, the residual variance of $y_t$ should increase. Mathematically, we prove that 
\begin{equation}
\det(\Sigma) = \det (\tilde{\Sigma}) \;. \label{eq:equaldets}
\end{equation} 

The proof of \eqref{eq:equaldets} is the only part that requires the analytic description of $\tilde{\Sigma}$, and is the main difficulty of the overall proof. It is not straightforward, because $\tilde{\Sigma}$ depends on the inverse of the covariance matrix $C_\bz (0)$, while we only have an analytic description of $\vech C_\bz (0)$. 
From \eqref{eq:equaldets}, it is easy to infer \mbox{$\tilde{\Sigma}_{yy} \leq \Sigma_{yy}$}, which completes the proof. 
{It is only in this final step that we need assumption (A2) that $\Sigma$ is diagonal.} 

\vspace{.5em}
\begin{IEEEproof}[Proof (Part 1: Proof that $\tilde{\Sigma}_{xx} \leq \Sigma_{xx}$)]
\label{sec:proofoftheoremVAR1}

As $A_1, \ldots, A_p$ are lower triangular matrices (assumption (A1)), $x_t$ is an autoregressive process of order $p$, 
\begin{equation}
x_t = a_1 x_{t-1} + \ldots + a_p x_{t-p} + \xi_t^x   \; \mbox{ with }  \; \Var( \xi_t^x)  = \Sigma_{xx} 
\end{equation}
Its time-reversed representation (cf. Section~\ref{sec:timereversal}) is
\begin{equation}
x_t = a_1 x_{t+1} + \ldots + a_p x_{t+p} + \tilde{\xi}_t^x , \; \mbox{ with }  \Var( \tilde{\xi}_t^x )  = \Sigma_{xx} 
\end{equation}

Because the unrestricted (or full) model \eqref{eq:VARprev} extends the restricted model by including $y_t$, \eqref{eq:abschaetzung1} follows:
\begin{equation}
\tilde{\Sigma}_{xx} \leq  \Var( \tilde{\xi}_t^x ) = \Sigma_{xx} \;.
\end{equation} 
\end{IEEEproof} 

\vspace{.5em}
\begin{IEEEproof}[Proof (Part 2: Proof that $\tilde{\Sigma}_{yy} \leq \Sigma_{yy}$)]  

As mentioned in the proof sketch, we need to derive \eqref{eq:equaldets}, the equality of the determinants \mbox{$\det \Sigma$} and \mbox{$\det \tilde{\Sigma}$} . To improve readability, we here treat only the case $p = 1$, and derive \eqref{eq:equaldets} for general $p \in \N \setminus \{ 0\}$ in Appendix~ \ref{app:proofequaldets}.

The proof relies on Sylvester's determinant theorem \cite{Akritas1996}, which states that for any matrices $K \in \R^{n \times m}$, $L \in \R^{m \times n}$:
\begin{equation}
\det (I + K L) = \det (I + L K) \;.
\label{eq:Sylvester}
\end{equation}
We then have: 
\begin{eqnarray*}
\det \Sigma &\overset{\eqref{eq:YuleWalker1}}{=} &\det (C_\bz(0)  - A_1 \cdot C_\bz(0)  \cdot A_1^\top) \\
&=  &\det (C_\bz(0)) \cdot \det(I   - A_1 \cdot C_\bz(0)  \cdot A_1^\top \cdot C_\bz(0)^{-1}) \\
&\overset{\eqref{eq:Sylvester}}{=}  &\det (C_\bz(0)) \cdot \det(I   -  C_\bz(0)  \cdot A_1^\top \cdot C_\bz(0)^{-1} A_1 ) \\
&= &\det (C_\bz(0)  - C_\bz(0)  \cdot A_1^\top \cdot C_\bz(0)^{-1} A_1  \cdot C_\bz(0) ) \\
&\overset{\eqref{eq:revA1E}}{=} &\det \tilde{\Sigma} \;.
\end{eqnarray*}
From the result of \emph{Part 1}~\eqref{eq:abschaetzung1}, the equality of residual covariance determinants~\eqref{eq:equaldets} (derived for general $p$ in Appendix~\ref{app:proofequaldets}), and assumption (A2) of uncorrelated residuals in $\Sigma$, we then obtain:
\begin{eqnarray}
\Sigma_{xx} \Sigma_{yy} \overset{(A2)}{=} \det{\Sigma}  \overset{\eqref{eq:equaldets}}{=}  \det{\tilde{\Sigma}} &=& \tilde{\Sigma}_{xx} \tilde{\Sigma}_{yy} - \tilde{\Sigma}_{xy} \tilde{\Sigma}_{xy} \nonumber \\
&\leq & \tilde{\Sigma}_{xx} \tilde{\Sigma}_{yy}  \nonumber \\
&\overset{\eqref{eq:abschaetzung1}}{\leq}& \Sigma_{xx}  \tilde{\Sigma}_{yy}  \nonumber \\
\Leftrightarrow   \Sigma_{yy}  &\leq&  \tilde{\Sigma}_{yy} \;. \nonumber
\end{eqnarray} 
\end{IEEEproof}

\emph{Strict inequality.} Let us further note that inequality \eqref{eq:Fnetdiff} for difference-based TRGC is strict in the presence of causal interaction. The following theorem holds.

\begin{theorem}
Under assumptions (A1)-(A3), it holds that
\begin{equation} 
\tilde{D}_{x \rightarrow y}^{(net)} = 0 \Leftrightarrow A_1, \ldots , A_p \mbox{ are diagonal}\, .
\end{equation}
\end{theorem}

The proof is provided in Appendix \ref{app:thm1InequalityDefiniteness}.  Combined with Corollary~1, Theorem~2 immediately implies that  $\tilde{D}_{x \rightarrow y}^{(net)} > 0$ in the presence of unidirectional information flow from $x_t$ to~$y_t$. That is, net Granger-causal information flow from $x_t$ to $y_t$ is truly reduced and cannot remain the same when the signal is time-reversed.

\section{Experiments}
\label{sec:experiments}

In this section, we provide an empirical investigation of model violations and other factors influencing the performance of Granger causal measures using numerical simulations. After describing the tested methods and performance measures (Section~\ref{sec:exp_methods}), we compare several variants of TRGC in either the presence or absence of noise (Section~\ref{sec:expvariants}). We then investigate the influence of common drivers, various types of noise (Section~\ref{sec:GCproblemsnoise} and \ref{sec:GCproblemsnoise_interaction}) and downsampling (Section~\ref{sec:GCproblemsdownsampling}) on standard Granger causality and Diff-TRGC. 

\subsection{Experimental setup}
\label{sec:exp_methods}

We consider bivariate time series in the presence of unidirectional information flow ($x_t \rightarrow y_t$) as well as in the absence of causal interaction.  Unless otherwise stated, time series of length $T=2000$ are generated from stationary VAR(5) processes, whose autoregressive coefficients are drawn from a normal distribution with mean~$0$ and standard deviation $\sigma_A = 0.2$. The absence of causal interaction is modeled by setting respective AR coefficients to zero. Residuals are generated from a normal distribution with diagonal covariance matrix, whose entries are drawn from the standard uniform distribution. 

We compare standard GC as well as Net-GC to Diff-TRGC (see~\eqref{eq:diffTRGC}). In Section~\ref{sec:expvariants}, we also include Conj-TRGC (see~\eqref{eq:conjTRGC}), the conjunction of Net-GC and Diff-TRGC (see~\eqref{eq:netdiffTRGC}), and  a variation of Diff-TRGC, in which $\tilde{D}_{x \rightarrow y}^{(net)}$ is computed using only the full bivariate models according to \eqref{eq:D_fullmodel}. This variant is denoted by \emph{Diff-TRGC (full)}. 

All statistical tests are performed at significance level \mbox{$\alpha = 0.05$}. For standard GC, we perform two separate F-tests, one to assess whether $x_t$ Granger-causes $y_t$, and one to assess whether $y_t$ Granger-causes $x_t$. It is possible that both variables are estimated to Granger-cause each other. In contrast, all other metrics indicate net directionality. We assess their statistical significance by bootstrapping residuals from the regression model: We regress $\bz_t$ on its past and future values $\bz_{t-p}, \ldots, \bz_{t-1}, \bz_{t+1}, \ldots \bz_{t+p}$, and retain the fitted values $\hat{\bz}_t$ and residuals $\hat{\epsilon}_t := \bz_t - \hat{\bz}_t$. In each bootstrap repetition, causality metrics are computed on synthetic variables $\bz_t^* := \hat{\bz}_t + \hat{\epsilon}_s$, where $s$ is selected randomly for each $t$. Percentile confidence intervals are then constructed from the bootstrap sampling distribution. Significance is determined by evaluating if the confidence interval does not contain 0. We use 500 bootstrap samples and select the number of lags $p$ as the optimizer of Schwarz's Bayesian Information Criterion (BIC) \cite{Schwarz1978}. 

All experiments are repeated 300 times. In each run, a true positive (TP) is defined as a significant detection of the true direction of interaction. The \emph{true positive rate} (TPR) is the fraction of true positives among all runs. 
It is here also referred to as the \emph{sensitivity} or \emph{power}. A false positive (FP) is defined as a significant detection of the wrong direction of interaction, or a significant detection of causal interaction in the absence of any causal interaction.  The \emph{false positive rate} (FPR) is the fraction of false positives among all tested runs.

\subsection{Comparison of TRGC variants under interaction}
\label{sec:expvariants}

\begin{figure}
 \centering
\subfigure[No measurement noise]{
 \includegraphics[width =0.49\linewidth]{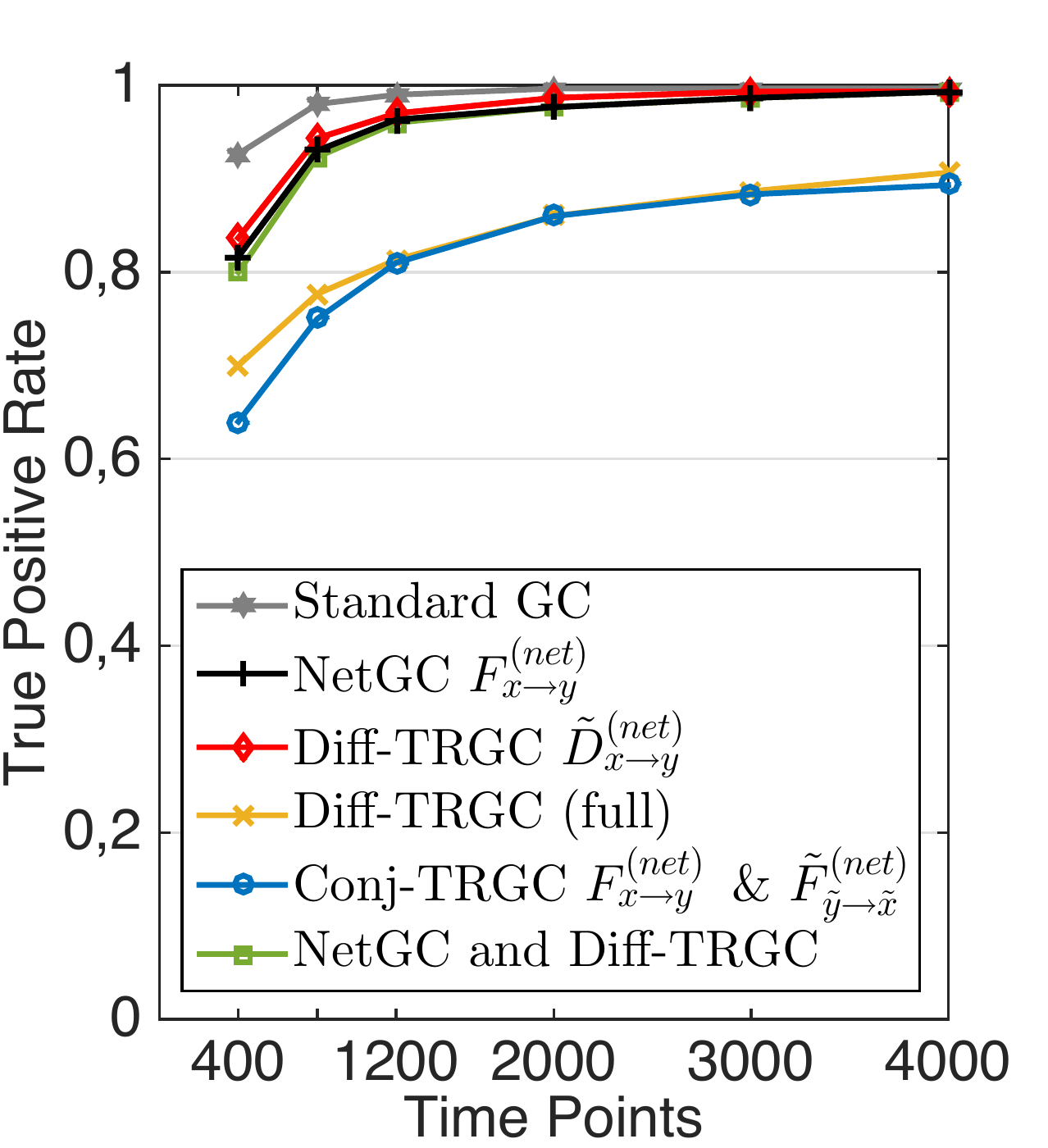} 
  \includegraphics[width =0.49\linewidth]{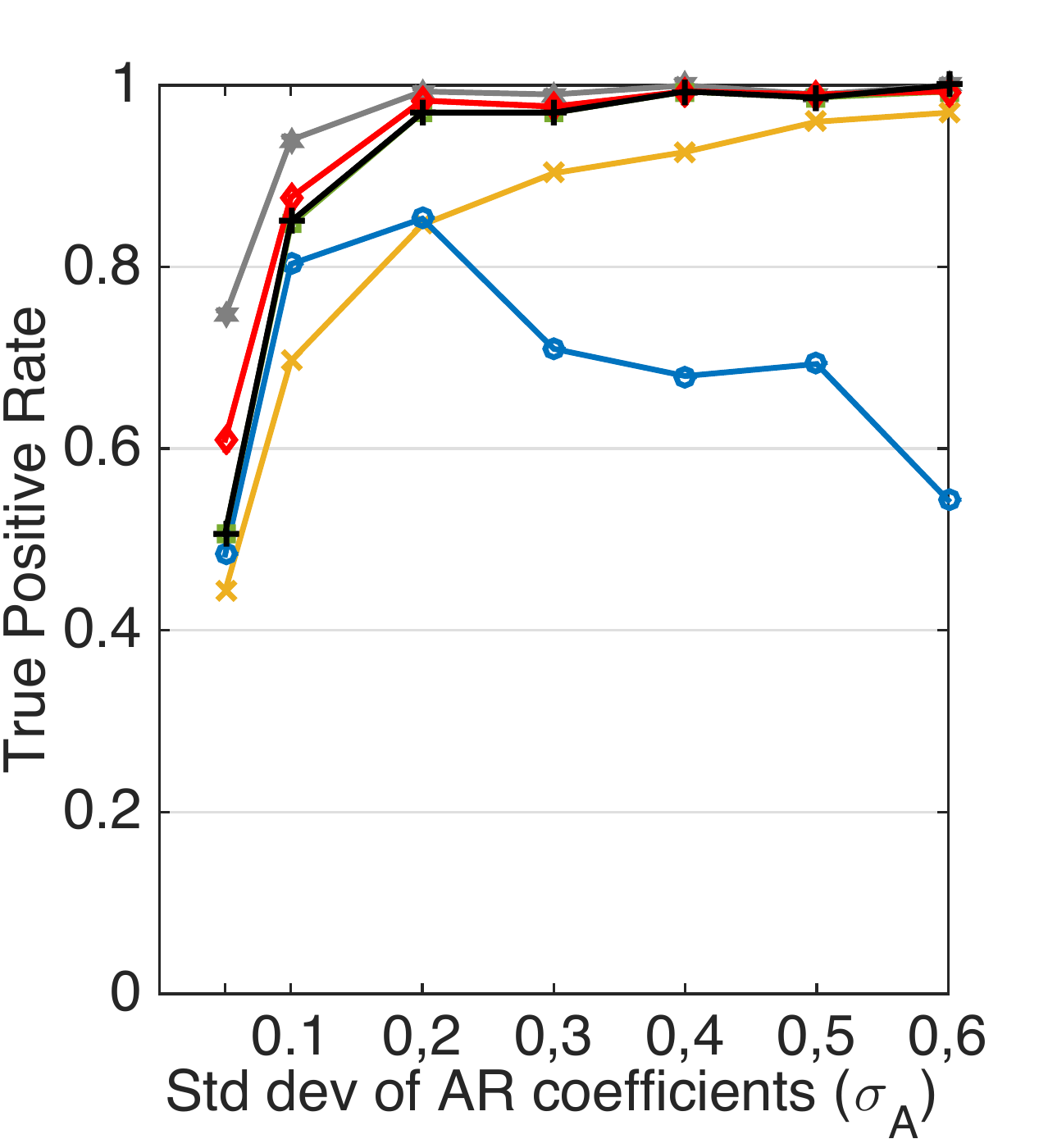}
}
\subfigure[Linearly mixed, auto-correlated noise]{
\includegraphics[width =0.98\linewidth]{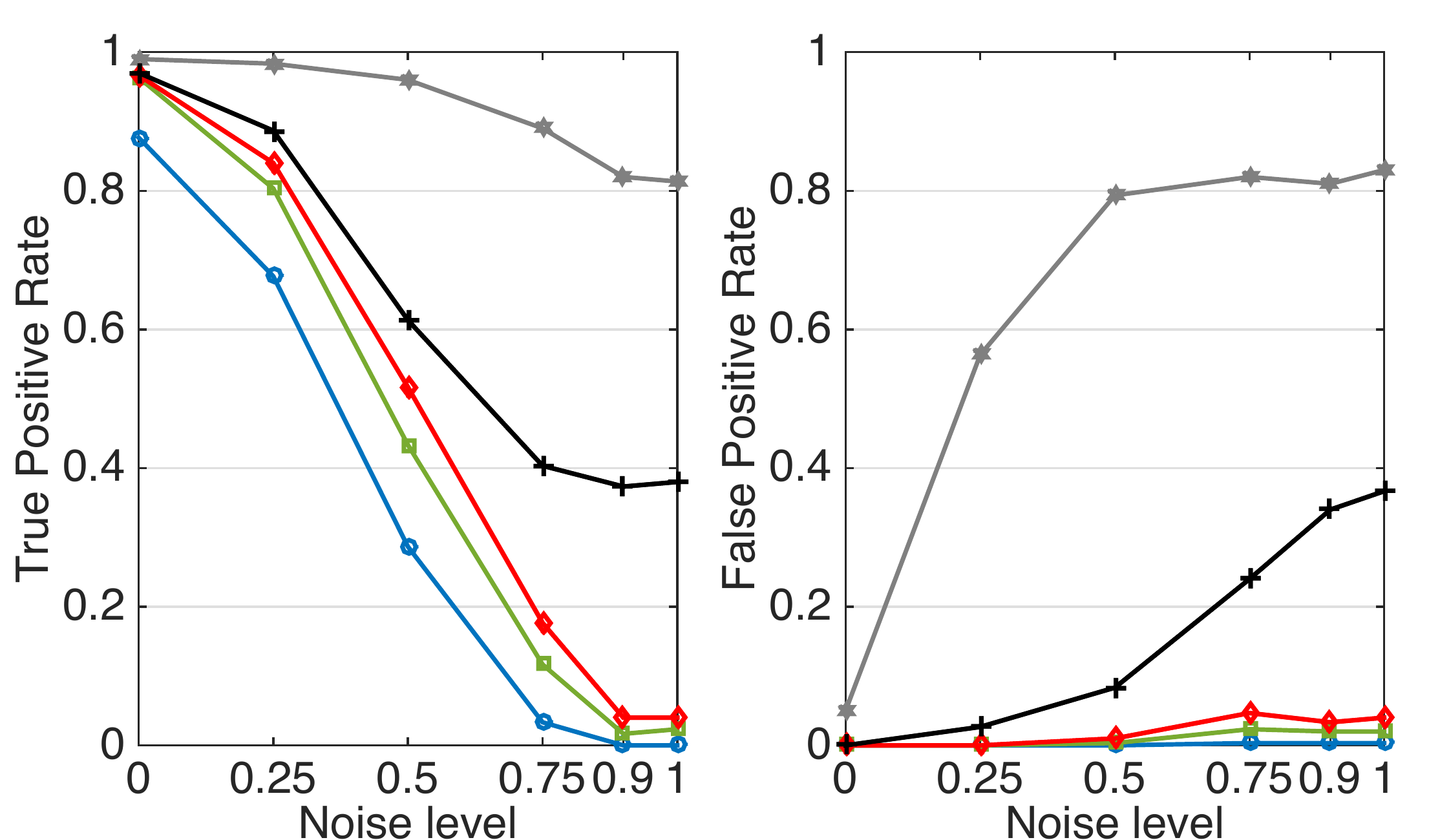} }
\caption{Performance of Granger causality and different variants of time-reversed Granger causality (TRGC). (a) True positive rate in the noiseless case as a function of the number of samples $T$ for fixed standard deviation $\sigma_A = 0.2$ of the AR coefficients, and as a function of $\sigma_A$ for fixed $T= 2000$. (b) True and false positive rates as a function of the SNR for additive mixed autocorrelated noise (according to~\eqref{eq:superposednoise}) for $T= 2000$ and $\sigma_A = 0.2$.}
\label{fig:simul_variants}
\end{figure} 

We assess Granger causality and time-reversed Granger causality in the presence of unidirectional interaction considering differing sample sizes, standard deviations of the AR parameters, noise types and signal-to-noise ratios (SNR). 

In a first experiment, we consider the noiseless case, and vary the sample size from 400 to 4000 for a fixed standard deviation $\sigma_A = 0.2$ of the AR coefficients. In a second experiment, we vary the standard deviation $\sigma_A$  at a constant sample size of $T = 2000$. This experiment thus tests the impact of the strength of the causal connections relative to the innovation noise. The standard deviations tested are 0.05, 0.1, 0.2, ..., and 0.6. Finally, for a fixed standard deviation $\sigma_A = 0.2$, and a fixed sample size $T=2000$, we add linearly mixed, autocorrelated measurement noise $\eta_t \in \R^2$ to each system according to
\begin{equation}
\begin{bmatrix} x_t \\ y_t \end{bmatrix} = (1- \gamma) \begin{bmatrix} x^{(l)}_t \\ y^{(l)}_t \end{bmatrix} + \gamma \cdot \eta_t \;,
\label{eq:superposednoise}
\end{equation}
where the subscript $^{(l)}$ denotes the underlying \emph{latent} variables and $\gamma$ defines the signal-to-noise ratio (SNR).  Noise $\eta_t$ is generated by multiplying two independent AR(5) time-series with a random matrix $B$, with $\det (B) = 1$. We consider the signal-to-noise ratios 0, 0.25, 0.5, 0.75, 0.9 and 1.

The TP and FP rates attained in the three experiments are depicted in Figure~\ref{fig:simul_variants}. From Figure~\ref{fig:simul_variants}(a), we see that \mbox{Diff-TRGC (full)}, which computes the difference score $\tilde{D}_{x \rightarrow y}^{(net)}$ only using the full model according to~\eqref{eq:D_fullmodel},  seems to be suboptimal for finite samples. While we have demonstrated the equivalence of \eqref{eq:D_fullmodel} to the original definition \eqref{eq:D_def} for infinite samples in Section~\ref{sec:TRGCdescription}, this equivalence does not hold for the finite samples studied here. Estimating residuals from the restricted models increased the power of the test for all investigated parameter settings. 

Conj-TRGC has lower power relative to Diff-TRGC. This is particularly so for high $\sigma_A$, which corresponds to a dominance of the dynamical and causal aspects of the model comprised in the AR coefficients relative to the innovation noise. This result is not unexpected, as time-reversing the signals does not necessarily reverse the direction of information flow. Note that, on the other hand, Conj-TRGC is the more conservative measure compared to Diff-TRGC and could be expected to produce fewer spurious results in the presence of noise. However, as we see in Figure~\ref{fig:simul_variants}(b), both variants yield almost no spurious results in the presence of measurement noise. We will therefore use Diff-TRGC in the remaining experiments.

\subsection{Impact of latent variables and measurement noise in the absence of causal interaction}
\label{sec:GCproblemsnoise}

\begin{figure*}
\centering
\bf  \hspace{2cm} A \hspace{6.5cm} B \hfill {~}\\ 
\includegraphics[width =0.4 \textwidth]{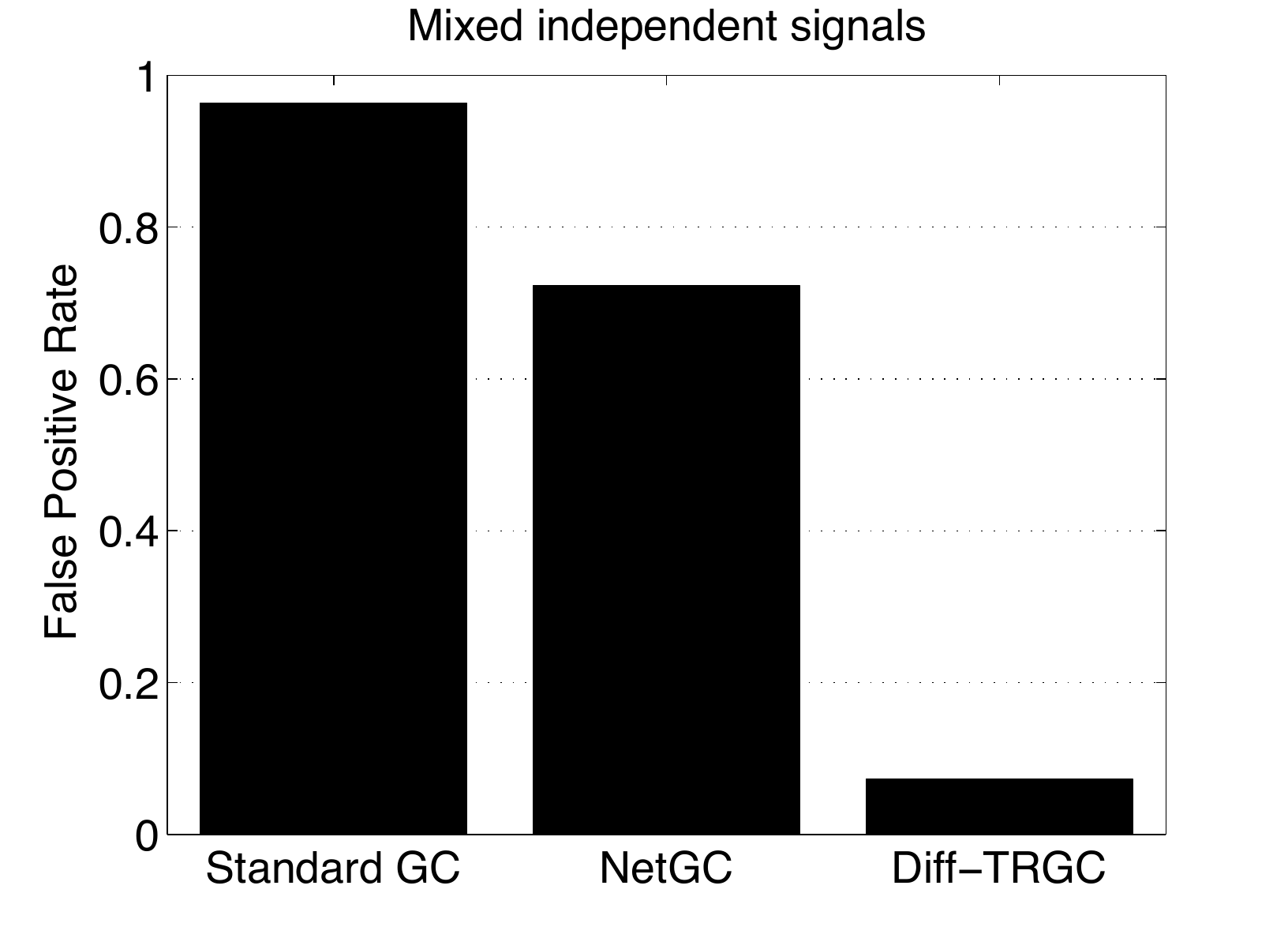}
 \includegraphics[width =0.4 \textwidth]{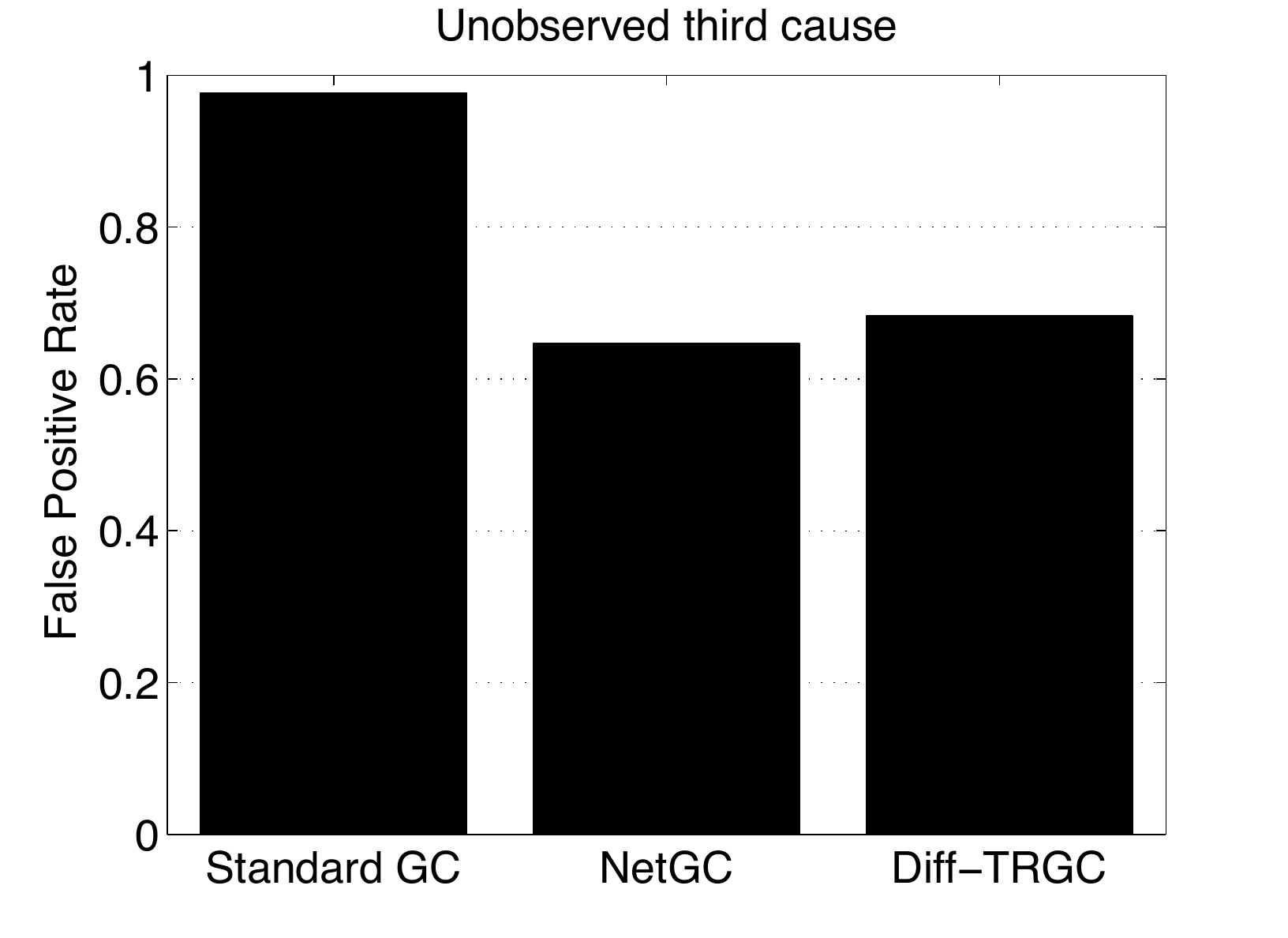}
\\ \hspace{2cm} C \hfill {~}\\ 
\includegraphics[width =0.91 \textwidth]{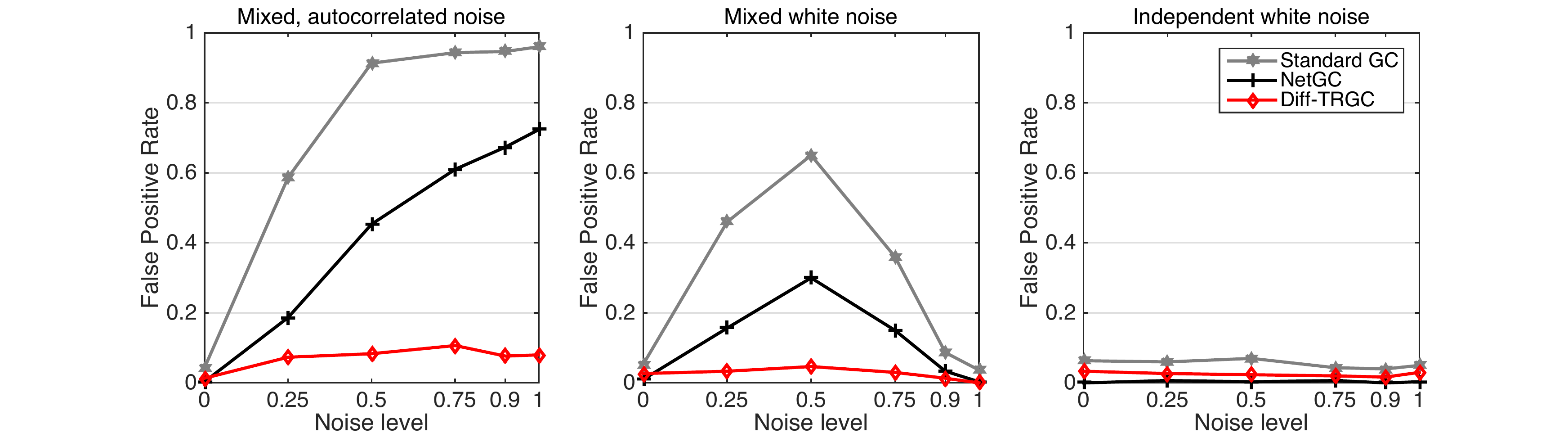}
\caption{False positive rates of Granger causality (standard GC and Net-GC) and difference-based time-reversed Granger causality (Diff-TRGC) as a function of the SNR for two signals lacking any causal connection.
\textbf{(A)} Instantaneous linear mixture of two independent univariate AR(5) processes. \textbf{(B)} Common unobserved cause. $x_t$ and $y_t$. 
\textbf{(C)} Superposition of two independent univariate AR(5) processes with additive Gaussian noise. 
}
\label{fig:failures_nocausality}      
\end{figure*}

Already Granger pointed out that standard Granger causality can lead to spurious results if not all relevant variables are incorporated in the model \cite{Granger1969}. In a bivariate system, GC cannot determine whether the observed variables $x_t$ and $y_t$ are both driven by a third common cause. This argument extends to multivariate systems, if a relevant confounding variable is not part of the measurement. 
Furthermore, standard GC is susceptible to \emph{measurement noise} \cite{Newbold1978,Nalatore2007,Nolte2008,Nolte2010,Sommerlade2012,Vinck2015} and to  \emph{instantaneous linear mixing} of activity, which is a major problem for example in the analysis of electroencephalographic (EEG) recordings \cite{GomezHerrero2008,Schoffelen2009,Haufe2013}. 
We demonstrate these effects here in additional simulations, in all of which no actual interaction occurs.
We consider three different scenarios.

\emph{(A) Linear mixing.} The observed time series $x_t$ and $y_t$ are a linear mixture of two independent signals $x^{(l)}_t$,  $y^{(l)}_t$, that is 
\begin{equation}
\begin{bmatrix} x_t \\ y_t \end{bmatrix} =M \begin{bmatrix} x^{(l)}_t \\ y^{(l)}_t \end{bmatrix} \;,
 \end{equation}
where $M \in \R^{2 \times 2}$ denotes the mixing matrix. $x^{(l)}_t$ and~$y^{(l)}_t$ were generated as two independent univariate AR(5) processes. 

\emph{(B) Common hidden cause.}
The observed time series $x_t$ and~$y_t$ are driven by a common unobserved cause $g_t$. Time series $x_t, y_t$, and $g_t$ are generated from a three-dimensional VAR(5) model with $\sigma_A = 0.3$, in which $g_t$ Granger-causes~$x_t$ and $y_t$, with no causal interaction between $x_t$ and $y_t$ as modeled by the respective AR coefficients being set to zero.  

\emph{(C) Additive noise.} 
The observed time series $x_t$ and $y_t$ are a superposition of two independent univariate AR(5) processes $x^{(l)}_t$,  $y^{(l)}_t$ and additive noise $\eta_t$ as in \eqref{eq:superposednoise}, with $\gamma \in \{0, 0.25, 0.5, 0.75, 0.9, 1 \}$ adjusting the SNR. 
We consider three different types of noise. \emph{Independent white noise} is generated from a normal distribution with diagonal covariance matrix, whose entries are drawn from the standard uniform distribution.  \emph{Mixed white noise} is created by multiplying independent noise with a random matrix $B$ with $\det (B) = 1$.  \emph{Mixed autocorrelated noise} is created by multiplying two independent AR(5) time-series with $B$.

Figure~\ref{fig:failures_nocausality} illustrates the behavior of standard Granger causality, Net-GC and Diff-TRGC in the various simulation settings. Values on the y-axis indicate the FP rate at significance level $\alpha = 0.05$. As all experiments are characterized by the absence of any interaction between $x_t$ and $y_t$, any significant detection of information flow either from $x_t$ to $y_t$ or $y_t$ to $x_t$ is counted as a false positive. 

It is apparent from Figure~\ref{fig:failures_nocausality} that standard GC and Net-GC lead to spurious detection of causality in all tested scenarios. Their behavior in the presence of noise (panel C) depends on the properties of that noise. Mixed noise (left and center plots of panel C) is generally very problematic, especially if it is also autocorrelated (left part). As $x_t$ and $y_t$ are already independent, adding independent noise (obviously) does not pose a problem here (right part of panel C).

In contrast to standard GC and Net-GC, time-reversed Granger causality implemented through Diff-TRGC is insensitive to mixtures of independent sources regardless of their spatial and temporal correlation structure (see panels A and C). This behavior thus reflects its known theoretical properties discussed in Section~\ref{sec:NoiseRobustness}. The presence of a hidden common confounder, however, cannot be ruled out by using time-reversed Granger causality (panel B).

\subsection{Impact of noise in the presence of causal interaction}
\label{sec:GCproblemsnoise_interaction}

\begin{figure*}
\centering
\bf \hspace{1cm} A \hspace{5cm} B \hspace{5cm} C \hfill \textcolor{white}{.}\\ 
\includegraphics[width =0.3\textwidth]{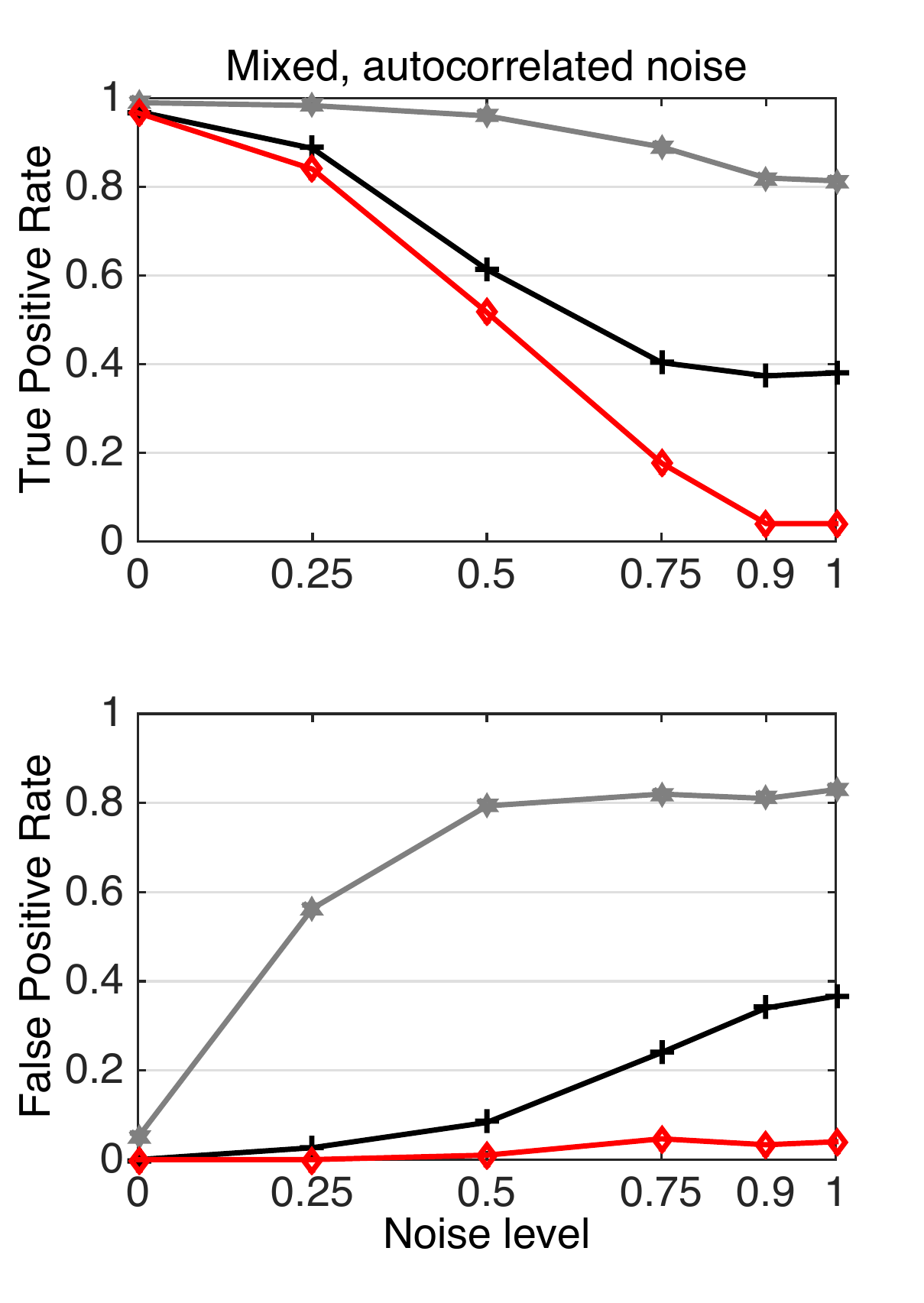} 
\includegraphics[width =0.3\textwidth]{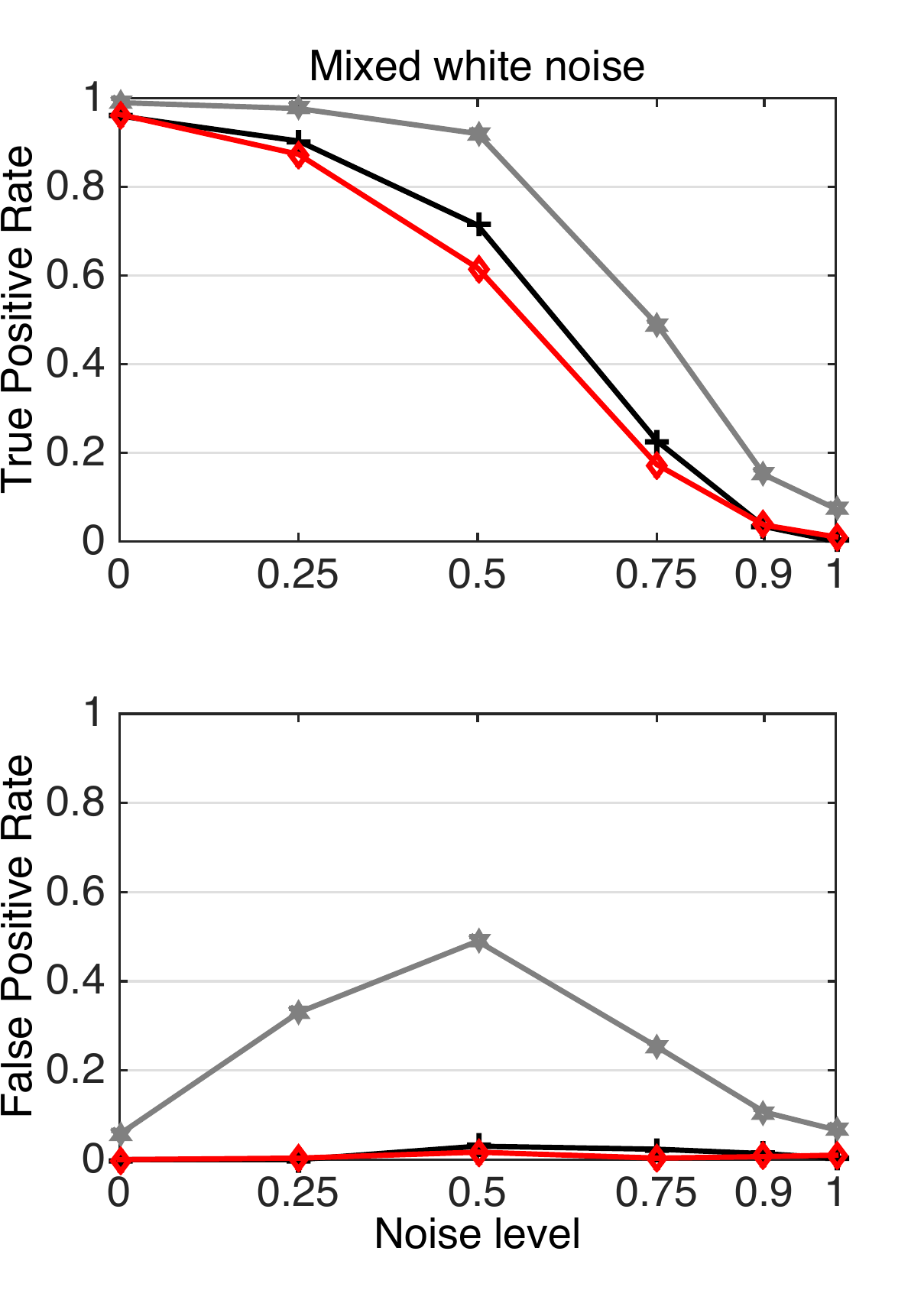}  
\includegraphics[width =0.3\textwidth]{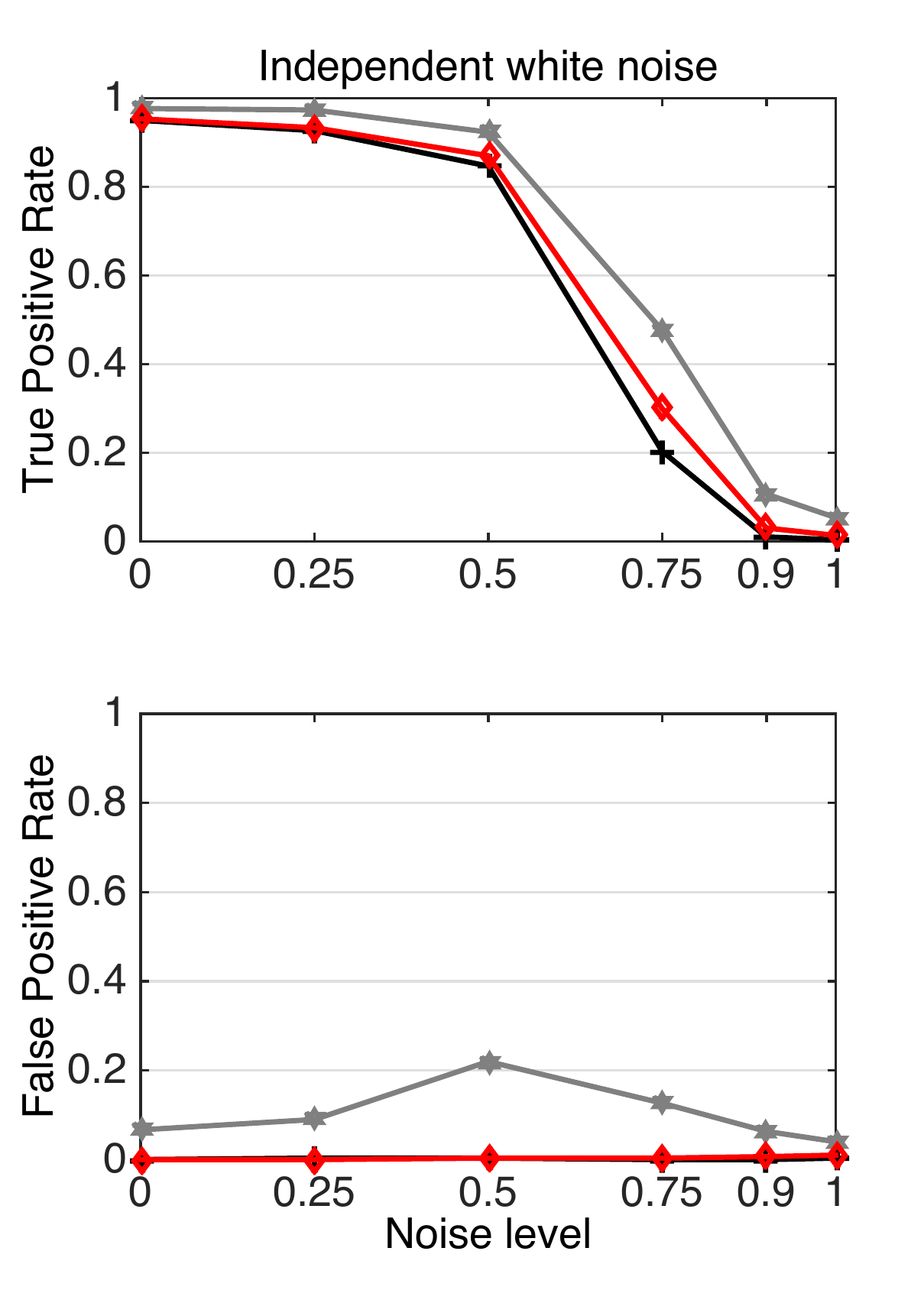}  \\  
\bf \hspace{1cm} D \hspace{5cm} E \hspace{5cm} F \hfill \textcolor{white}{.}\\ 
\includegraphics[width =0.3\textwidth]{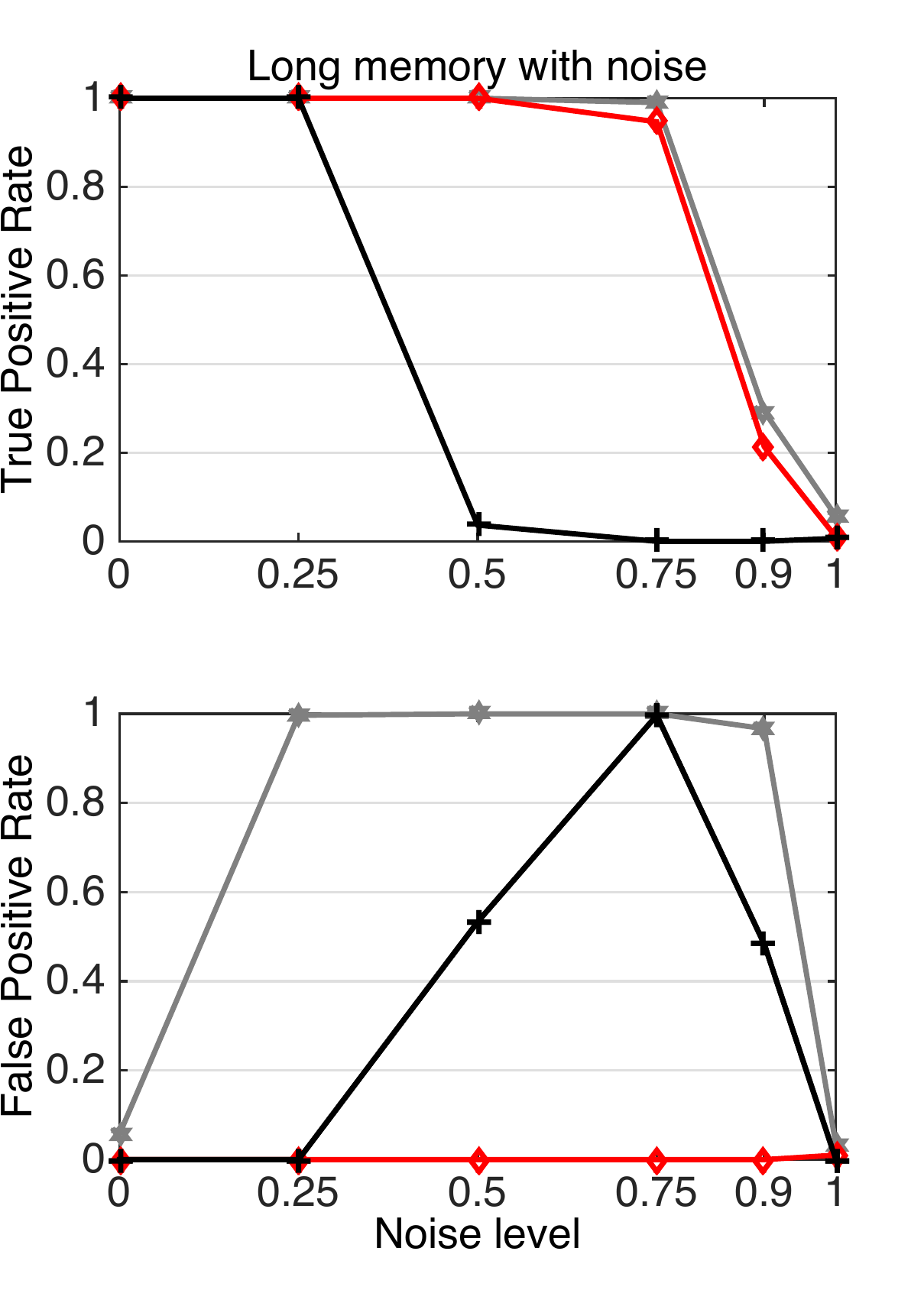}
\includegraphics[width =0.3\textwidth]{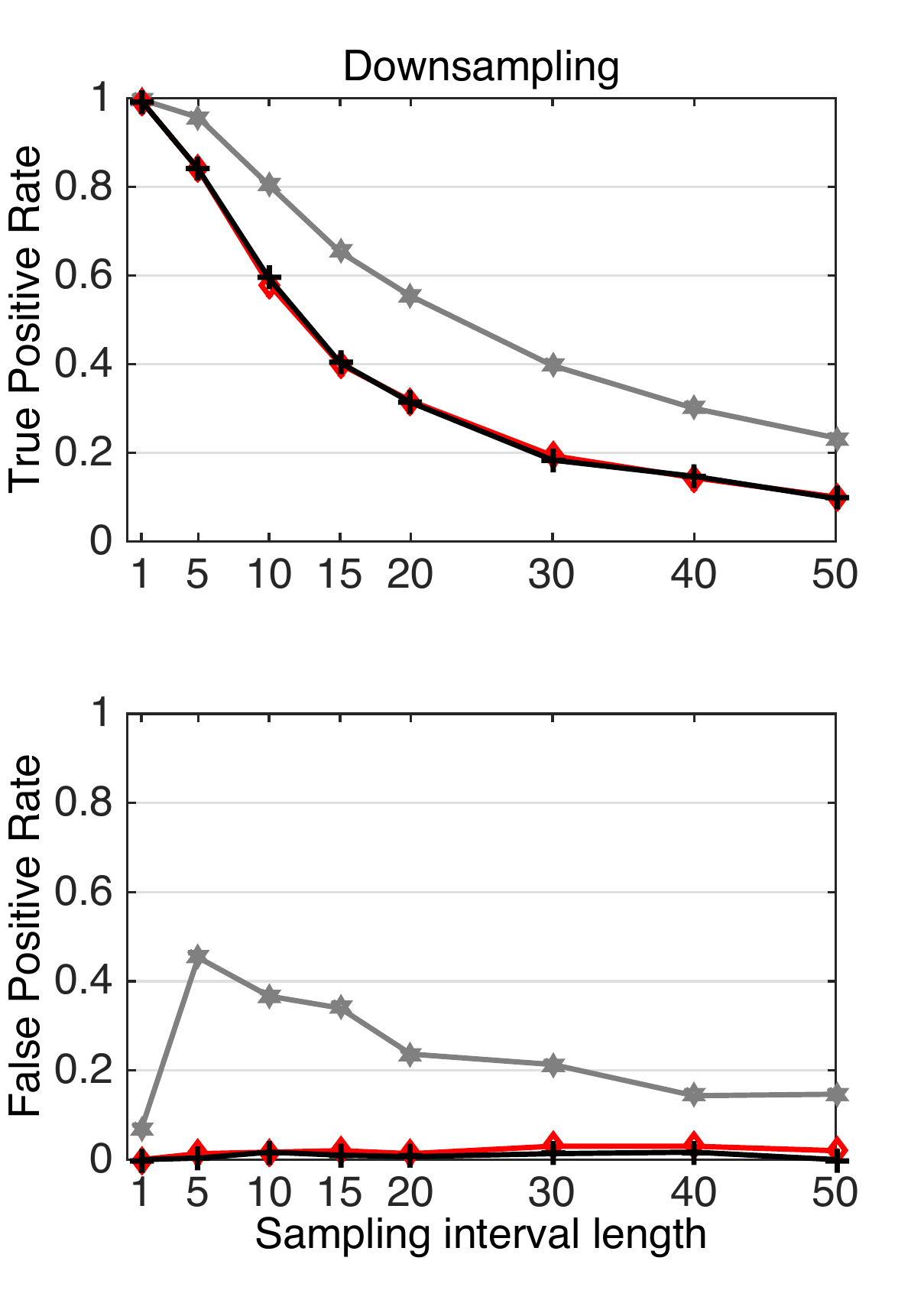}
\includegraphics[width =0.3\textwidth]{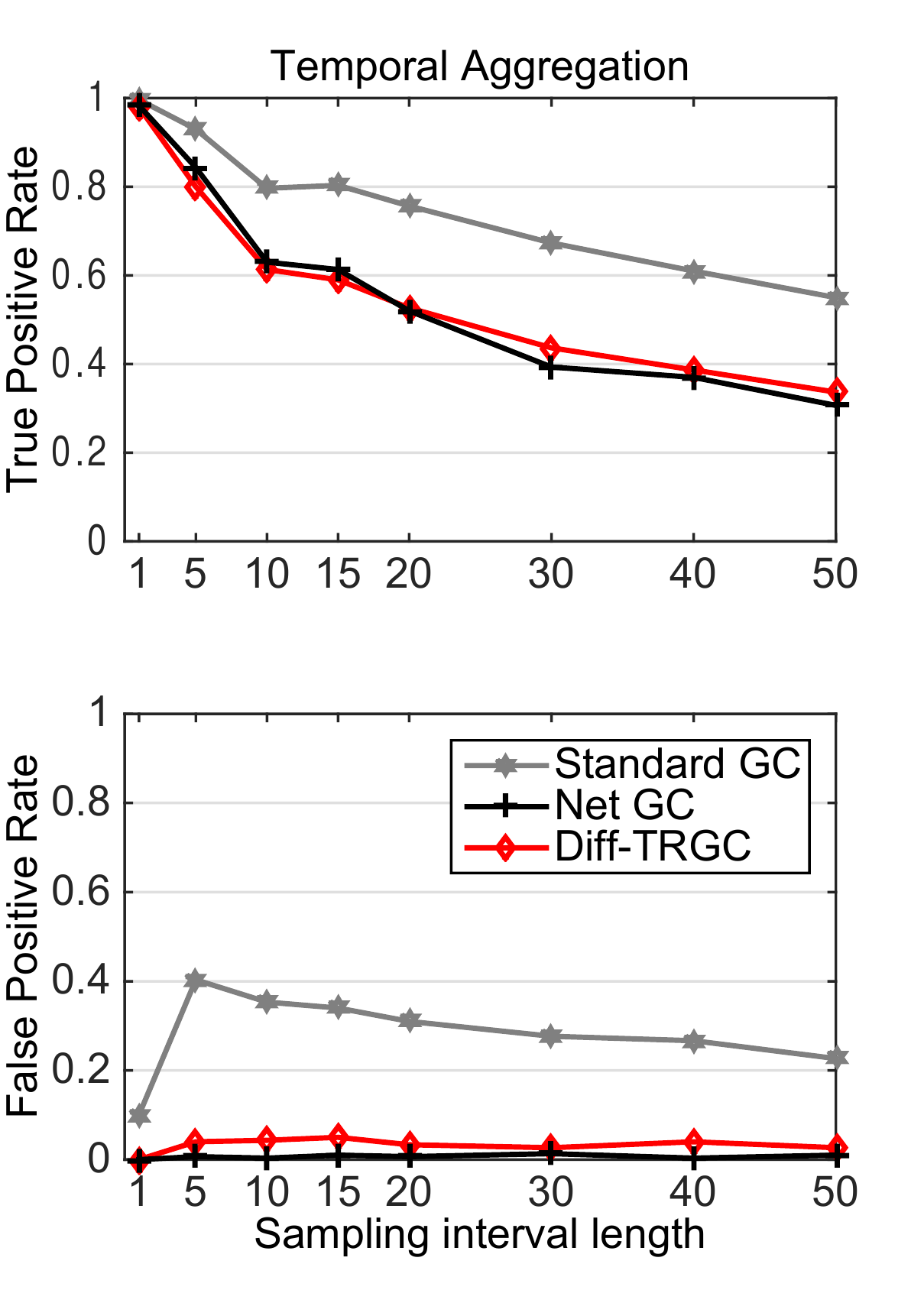} 
\caption{Performance of Granger causality (standard GC and Net-GC) and difference-based time-reversed Granger causality (Diff-TRGC) for two signals with unidirectional information flow from $x_t$ to $y_t$. Shown are the fractions of true positives ($x_t \rightarrow y$ detected) and false positives ($y_t \rightarrow x_t$ detected), when $x_t$ and $y_t$ are corrupted by noise (A-D), downsampling (E), and temporal aggregation (F). The underlying latent signals $x^{(l)}$ and $y^{(l)}$ were generated from VAR(5) processes with random AR coefficients, except for D, in which signals follow a VAR(1) process with long memory according to \eqref{eq:NolteVAR1}. 
}
\label{fig:failures_summary}
\end{figure*}

We further study the behavior of standard GC, Net-GC and Diff-TRGC in the presence of unidirectional causal interactions superimposed with noise. Four different scenarios are considered. In all cases, data are generated according to~\eqref{eq:superposednoise} with $x_t^{(l)}$ Granger-causing $y_t^{(l)}$. In the first three scenarios, (A-C), interacting signals from bivariate VAR(5) models are superimposed with noise. As in 
Section~\ref{sec:GCproblemsnoise}, we use mixed autocorrelated noise (scenario A), mixed white noise (B), and independent white noise (C). The same signal to noise ratios as in Section~\ref{sec:GCproblemsnoise} are used.

In the fourth scenario, (D), we simulate the following VAR(1) process with long memory:
\begin{equation}
\begin{aligned}  \begin{bmatrix} x_t^{(l)} \\ y_t^{(l)}  \end{bmatrix} &=  \begin{bmatrix} 0.95 & 0 \\ 1 & 0.5 \end{bmatrix} \begin{bmatrix} x_{t-1}^{(l)} \\ y_{t-1}^{(l)}  \end{bmatrix}  + \epsilon_t & \epsilon_t \sim \mathcal{N}(0, I) \\
x_t & = (1 - \gamma) \cdot x_t^{(l)} + \gamma \cdot \eta_t &  \eta_t \sim \mathcal{N}(0, 1)\\
y_t & = y_t^{(l)} \;, & 
\end{aligned}
\label{eq:NolteVAR1}
\end{equation} 
adopted from \cite{Nolte2010}, where $\mathcal{N}$ denotes the normal distribution. 

True positive and false positive rates as estimated from 300 simulation runs are reported in Figure~\ref{fig:failures_summary} (A-D). Just as in the absence of causality (cf. Section~\ref{sec:GCproblemsnoise}), we observe that linearly mixed, autocorrelated noise leads to the highest numbers of false detections for standard GC, while independent white noise leads to lowest FP rates. Diff-TRGC is characterized by negligible amounts of false positives in all cases at the cost of slightly decreased sensitivity as compared to standard GC in scenarios (A-C). Interestingly, Net-GC behaves very similar to Diff-TRGC in the presence of non-autocorrelated noise both in terms of sensitivity and specificity (B-C). In these settings, spurious causality could already be almost entirely eliminated  by testing for net Granger causality. This result, however, does not imply that Net-GC cannot be affected by non-autocorrelated noise in general. A counterexample is the system with long memory studied in scenario (D). Here, Net-GC (as well as standard GC) fails, because $y_t$ contains delayed but cleaner information about $x_t^{(l)}$ than $x_t$ itself and thus may help to predict future $x_t$. Diff-TRGC, however, robustly identifies $x_t$ as the driver.

Our examples show time-reversed Granger causality almost completely eliminates spurious causalities arising from any kind of additive noise. At the same time, it exhibits similar statistical power as net Granger causality. We also observe that net Granger causality is typically more robust with respect to additive noise than standard Granger causality. 

\subsection{Impact of downsampling and temporal aggregation}
\label{sec:GCproblemsdownsampling}

Spurious Granger causality has also been reported to arise due to downsampling and temporal aggregation \cite{Tiao1976,McCrorie2006,Zhou2014}, posing serious problems, for example, in functional magnetic resonance imaging (fMRI) \cite{Seth2013,Smith2011}. 

We generate data using a VAR(5) model with random coefficients with $\sigma_A = 0.3$, in which $x_t$ Granger-causes $y_t$. These data are decimated at different factors $\tau$ in two ways. In the downsampling scenario (E), causal measures are applied to time series of length $T = 2000$ constructed from the original time series by skipping $\tau -1$ time points in between sampled data points. In the temporal aggregation scenario (F), time series of length $T = 2000$ are constructed from the original time series by averaging over $\tau$ data points. No noise was added. 

Figure~\ref{fig:failures_summary} (E-F) depicts TP and FP rates attained in the two scenarios as a function of $\tau$. We see that Net-GC and Diff-TRGC are more robust then standard GC. Both Net-GC and Diff-TRGC did not result in spurious causality.

\section{Discussion}
\label{sec:Discussion}

We established the theoretical guarantee that difference-based time-reversed Granger causality (Diff-TRGC) indicates the correct direction of causality in bivariate autoregressive processes characterized by unambiguous unidirectional information flow. Our results complement previous work by \cite{Haufe2013,Haufe2012} showing that TRGC in general correctly rejects causal interpretations for mixtures of non-interacting sources (thus, in the absence of any causality). While further compelling intuitive ideas for robust causality measures have been presented \cite{Nolte2008,Vicente2011,Vinck2015}, our result provides, to the best of our knowledge, the first proof of the correctness of one of such techniques (Diff-TRGC) for a relatively general class of time-series models. 

Our theory is accompanied by simulations, in which we confirmed that time-reversed Granger causality robustly detects the presence of true causal interactions in various realistic scenarios including mixed noise and downsampling. We showed that Diff-TRGC is able to infer correct directionality with similar power as net Granger causality, while at the same time producing fewer (in most cases, negligible amounts of) false alarms than Net-GC and standard GC. We therefore suggest to use Diff-TRGC whenever the data under study are likely to be corrupted by noise.

\subsection{Correlated residuals}
\label{sec:corrRes}

To define an unambiguous uni-directional information flow, our theory assumes uncorrelated residuals, as is common in the literature. Correlated residuals indicate instantaneous effects that the variables exert on each other. While we would not expect correlated residuals if the VAR model accurately describes the data generating process, such effects are likely to occur in practice (e.g., if the sampling rate of the acquired data falls below the time scale of the causal interactions). They pose severe problems for causal estimation, because they can be explained by several possible data generating models, the coefficients of which cannot be uniquely identified using second order information only. 

\emph{Data generating models.}
Instantaneous interactions can be modeled implicitly through correlated residuals in classical VAR processes, or explicitly, for example using so-called `structural' VAR (SVAR) processes \cite{Luedtkepohl2007,Hyvarinen2010,Moneta2011}. By augmenting the VAR model with an instantaneous mixing matrix $\Gamma_0$, the SVAR model
\begin{equation} 
\bz_t = \sum_{h=0}^p \Gamma_h \bz_{t-h} + \bar{\epsilon}_t \;,
\label{eq:svarmodel}
\end{equation}
achieves that the residuals $\bar{\epsilon}_t$ are uncorrelated. Here, the diagonal of $\Gamma_0$ is assumed to be zero.

Correlated residuals emerge naturally in electrophysiological neuroimaging data, where the signals observable at the sensors (e.g., EEG electrodes) are a linear mixture of the latent activity of possibly interacting neuronal populations within the brain. A model for such mixtures of potentially interacting sources is given by 
\begin{equation}
\begin{aligned} 
\bz_t & = M \bz_t^{(l)} \, , \quad 
\bz_t^{(l)} & = \sum_{h=1}^p B_h \bz_{t-h}^{(l)} + \acute{\epsilon}_t \;,
\end{aligned}
\label{eq:scsamodel}
\end{equation}
where $\bz_t\in \mathbb{R}^d$ denotes the observed data, $\bz_t^{(l)}  \in \mathbb{R}^d$ denotes the activity of underlying latent variables (e.g., brain sources) following a VAR($p$) process with uncorrelated residuals $\acute{\epsilon}_t$, and $M \in \mathbb{R}^{d \times d}$ is an unknown mixing matrix (representing, e.g., the volume conduction effect of the human head). We call \eqref{eq:scsamodel} the `mixture of interacting sources' model. 

Note that VAR models with correlated residuals, SVAR models, and mixture of interacting sources models can be used interchangeably to represent the same statistical process. For example, an interacting sources model~\eqref{eq:scsamodel} can be equivalently written as 
a VAR($p$) process~\eqref{eq:VARp} with coefficients
\begin{equation}
A_h = M B_h M^{-1} \quad , \quad h \in \{1, \hdots, p\}
\label{eq:scsamodel_var_params}
\end{equation}
and correlated residuals $\epsilon_t = M \tilde{\epsilon}_t$. Likewise, an SVAR($p$) process~\eqref{eq:svarmodel} can be converted into a VAR($p$) process with correlated residuals $\epsilon_t = (I-\Gamma_0)^{-1} \bar{\epsilon}_t$ and coefficients
\begin{equation}
A_h = (I-\Gamma_0)^{-1}  \Gamma_h \quad , \quad h \in \{1, \hdots, p\} \;.
\label{eq:svar_model_var_params}
\end{equation}
The reverse transformations from VAR models to SVAR or interacting source models, as well as the transformations between SVAR and interacting source models, are, however, not unique (see \emph{Model identifiability}).

Ambiguous causal interpretations can emerge in cases where one of the three models indicates time-delayed causal interactions through non-zero off-diagonal coefficients in the $A_h$, $B_h$ or $\Gamma_h$, while another one does not. This ambiguity can in general only be resolved if the model generating the data is known a-priori. In case of EEG data, for example, \eqref{eq:scsamodel} reflects the true data-generating process. Therefore, only the parameters $B_h$ of the source VAR process~\eqref{eq:scsamodel} permit meaningful causal interpretation (wrt. to the source variables~$\bz_t^{(l)} $), while, for example,
the VAR parameters in \eqref{eq:scsamodel_var_params} are distorted by the mixing matrix $M$.

\emph{Model identifiability.} 
A further complication in the presence of instantaneous effects in the data is that for mixture of interacting sources as well as SVAR models, the parameters are not uniquely defined from second order information only. 
This can be best seen for the latter model~\eqref{eq:scsamodel}. Identifying the model parameters requires the estimation of a full factorization of the data into a mixing matrix $M$ and source time series~$\bz_t^{(l)}$. 
This means that the estimation problem falls into the blind source separation~(BSS) setting, in which Gaussianity of the factors is not sufficient for their identification. The classical approach to BSS, independent component analysis~(ICA) assumes statistical independence and non-Gaussianity of the sources $\bz_t^{(l)}$ to ensure identifiability. This concept can be adopted in the context of source AR models by enforcing independence/non-Gaussianity of the residuals of the source AR process in \eqref{eq:scsamodel} \cite{GomezHerrero2008,Haufe2010,Chiang2012}. In a similar way, independence of residuals has been used in the identification of SVAR models \cite{Hyvarinen2010,Peters2013}.

\emph{Example.}
Consider the following VAR(1) process with correlated residuals: 
 
{\small \setlength{\arraycolsep}{2pt}
\begin{flalign*}
 \begin{bmatrix} x_t \\ y_t \end{bmatrix}  &=\begin{bmatrix} 0.7 & 0 \\ -0.12 & 0.9 \end{bmatrix}  \begin{bmatrix} x_{t-1} \\ y_{t-1} \end{bmatrix} + \epsilon_t  ,\quad \langle \epsilon_t \epsilon_t^\top \rangle = \begin{bmatrix} 1 & 0.6 \\ 0.6 & 1 \end{bmatrix}  \,. &
\end{flalign*} \setlength{\arraycolsep}{5pt}}

\noindent This process can also be represented by the SVAR(1) model

{\small \setlength{\arraycolsep}{2pt} \begin{flalign*}
 \begin{bmatrix} x_t \\ y_t \end{bmatrix}  &=   \begin{bmatrix} 0 & 0 \\ 0.6 & 0 \end{bmatrix}  \begin{bmatrix} x_{t} \\ y_{t} \end{bmatrix} +  \begin{bmatrix} 0.7 & 0 \\ -0.54 & 0.9 \end{bmatrix}  \begin{bmatrix} x_{t-1} \\ y_{t-1} \end{bmatrix} + \bar{\epsilon}_t &
\end{flalign*} \setlength{\arraycolsep}{5pt}} 

\noindent  as well as the mixtures of interacting sources model
 
{\small \begin{equation*}
\arraycolsep=2pt
\begin{aligned}
\begin{bmatrix} x_t \\ y_t \end{bmatrix} &= \begin{bmatrix} 1 & 0 \\ 0.6 & 0.8 \end{bmatrix} \begin{bmatrix} x_{t}^{(l)} \\ y_{t}^{(l)} \end{bmatrix} , \quad  \begin{bmatrix} x_{t}^{(l)} \\ y_{t}^{(l)} \end{bmatrix} = \begin{bmatrix} 0.7 & 0 \\ 0 & 0.9 \end{bmatrix}  \begin{bmatrix} x_{t-1}^{(l)} \\ y_{t-1}^{(l)} \end{bmatrix} + \acute{\epsilon}_t ,
\label{eq:scsaexamp}
\end{aligned}
\end{equation*}}

\noindent  with uncorrelated residuals {\small $\langle \bar{\epsilon}_t \bar{\epsilon}_t^\top \rangle  = \begin{bmatrix} 1 & 0 \\ 0 & 0.64 \end{bmatrix} $, $ \langle \acute{\epsilon}_t \acute{\epsilon}_t^\top \rangle = I$}. Note that both the VAR(1) and the SVAR(1) representation indicate unidirectional causal interaction between the observed variables $x_{t}$ and $y_{t}$, whereas the mixture model suggests that the observed data can also arise from a mixture of two independent latent sources $x_{t}^{(l)}$ and $y_{t}^{(l)}$. However, another equivalent mixture model 

{\small
\begin{equation*}
\arraycolsep=1pt
\begin{aligned}
\begin{bmatrix} x_t \\ y_t \end{bmatrix} &= \begin{bmatrix} -\sqrt{0.2} & \sqrt{0.8} \\ \sqrt{0.2} & \sqrt{0.8} \end{bmatrix} \begin{bmatrix} x_{t}^{(l)} \\ y_{t}^{(l)} \end{bmatrix} , \begin{bmatrix} x_{t}^{(l)} \\ y_{t}^{(l)} \end{bmatrix} = \begin{bmatrix} 0.86 & 0.08 \\ 0.08 & 0.74 \end{bmatrix}  \begin{bmatrix} x_{t-1}^{(l)} \\ y_{t-1}^{(l)} \end{bmatrix} + \tilde{\tilde{\epsilon}}_t 
\label{eq:scsaexamp2}
\end{aligned}
\end{equation*}}

\noindent  with $\langle \tilde{\tilde{\epsilon}}_t \tilde{\tilde{\epsilon}}_t^\top \rangle = I$ suggests bidirectional informational flow on the source level. Similarly, the following SVAR(1) model indicates bidirectional flow

\small
\setlength{\arraycolsep}{2pt}
\begin{flalign*}
 \begin{bmatrix} x_t \\ y_t \end{bmatrix}  &=  \begin{bmatrix} 0 & 0.6 \\ 0 & 0 \end{bmatrix}  \begin{bmatrix} x_{t} \\ y_{t} \end{bmatrix} +  \begin{bmatrix} 0.772 & -0.54 \\ -0.12 & 0.9 \end{bmatrix}  \begin{bmatrix} x_{t-1} \\ y_{t-1} \end{bmatrix} + \bar{\bar{\epsilon}}_t \quad,  &\\
\langle \bar{\bar{\epsilon}}_t \bar{\bar{\epsilon}}_t^\top \rangle & =  \begin{bmatrix} 0.64 & 0 \\ 0 & 1 \end{bmatrix} . &
\end{flalign*} 
\setlength{\arraycolsep}{5pt}
\normalsize

\subsection{Future work}
\label{sec:futurework}

Further effort is required to investigate the behavior of TRGC in the presence of bidirectional information flow. Also, our theoretical analysis only covers the bivariate framework. Both Granger causality and TRGC can result in spurious causality when relevant variables are not included (cf. Fig.~\ref{fig:failures_nocausality}, panel B). Therefore, an extension of the analysis of time-reversal to general multivariate signals would be very interesting. Furthermore, it would be desirable to obtain theoretical guarantees for the performance of TRGC in the presence of true interaction superimposed by noise in the form of bounds on the false positive rate. A major difficulty here is to obtain the residual covariance of the superposition of a VAR process and additive noise. {Analytically computing Granger causality in the presence of noise is mathematically involved even for special cases \cite{Nalatore2014}.} 

Finally, \cite{Haufe2013} showed that for any causality measure based on cross-covariances, differences of the scores obtained on the original and time-reversal signals correctly indicate the absence of causality on mixtures of independent sources. While we focused here on Granger causality, it remains to be shown whether validity of time-reversal in the presence of causal interaction can also be demonstrated for other causality measures.

\appendices
\small

\section{Proofs for VAR($p$)}
\label{app:proofs}

\subsection{The VAR($p$) process and its cross-covariance function}
\label{app:YuleWalker}

Consider a stable bivariate VAR($p$) process,  $\bz_t  \in \R^2$, as defined in~\eqref{eq:VARp}, 
\begin{equation*}
\bz_t = A_1 \bz_{t-1} + A_2 \bz_{t-2} + \ldots + A_p \bz_{t-p} + \epsilon_t \;, 
\end{equation*}
where $\epsilon_t \in \R^2$ is a $2$-dimensional white noise process (i.e. $\langle \epsilon_t \rangle = 0$,  $\langle \epsilon_t \epsilon_{t-h}^\top \rangle = 0$ for $h  \in \Z \setminus \{0\}$, and $\langle \epsilon_t \bz_{t-h}^\top \rangle = 0$ for $h \in \N \setminus \{0\}$ ) with residual covariance matrix $\Sigma = \langle \epsilon_t \epsilon_t^\top \rangle$. 

Many results on VAR($1$) processes can be extended to higher order VAR($p$) processes by considering their VAR(1) form. Given the $2$-dimensional VAR($p$) process $\bz_t$, the corresponding $2p$-dimensional VAR(1) representation is defined as
\begin{equation}
Z_t = \bold A Z_{t-1} + E_t \label{eq:VAR1representation} \;,
\end{equation}
with 
\begin{equation*} \arraycolsep=2pt
Z_t = \begin{bmatrix} \bz_t \\ \bz_{t-1} \\ \vdots \\ \bz_{t-p+1} \end{bmatrix} \;, \bold A = \begin{bmatrix} A_1 & A_2 &\cdots & A_{p-1} & A_p \\ I & 0 & \cdots & 0 & 0 \\  0& I & \cdots & 0 & 0  \\   &  & \ddots & &  \vdots  \\  0& 0& \cdots & I & 0\end{bmatrix} \;, E_t = \begin{bmatrix} \epsilon_t \\ 0\\ \vdots \\ 0 \end{bmatrix} \;,
\end{equation*}
and residual covariance matrix
\begin{equation}
\arraycolsep=2pt
\Sigma_E = \langle E_t E_t^\top \rangle = \begin{bmatrix} \Sigma & 0 &\cdots & 0 \\ 0 & 0 & \cdots & 0  \\ \vdots  &  & \ddots & \vdots  \\  0& 0& \cdots & 0\end{bmatrix} \;.
\end{equation}
The cross-covariances of $Z_t$ are linked to the cross-covariances of $\bz_t$ through
\begin{equation} 
\arraycolsep=2pt
\medmuskip=1mu
C_Z(h) = \begin{bmatrix} C_\bz(h) & C_\bz(h+1) &\cdots & C_\bz(h+p-1) \\ C_\bz(h-1) & C_\bz(h) & \cdots & C_\bz(h+p-2)  \\  \vdots &  \vdots & \ddots &  \vdots  \\ C_\bz(h-p+ 1)& C_\bz (h-p + 2)& \cdots & C_\bz(h) \end{bmatrix} 
\end{equation}
for all $h \in \Z$. The Yule-Walker equations can then be expressed as 
\begin{eqnarray} 
C_Z(0) &=& \bold A \cdot C_Z(0) \cdot \bold A^\top + \Sigma_E  \quad \text{and} \label{eq:YuleWalkerVARp} \\  
C_Z(h) &=& \bold A \cdot C_Z(h-1) \qquad \forall h \in \N \setminus \{0 \} \label{eq:YuleWalkerVARp_h} \;.
\end{eqnarray}

Given $A_1, \ldots, A_p$, and $\Sigma$, the cross-covariances are uniquely determined:  Equation~\eqref{eq:YuleWalkerVARp} implies that \mbox{$ \vech(C_Z(0)) = (I - \bold A \kron \bold A)^{-1} \vech{\Sigma_E} $}, while $C_Z(h)$ for $h > 1$ can be recursively computed using \eqref{eq:YuleWalkerVARp_h}. Conversely, $A_1$, \ldots, $A_p$ and $\Sigma$ are uniquely determined by the cross-covariances through 
\begin{eqnarray} 
[ A_1 , A_2 ,\cdots , A_p]  = [C_\bz(1) , C_\bz(2) , \cdots , C_\bz(p) ] \cdot C_Z(0)^{-1}  
\end{eqnarray}
and
\begin{eqnarray} 
\Sigma = C_\bz(0) - [ A_1 , A_2 ,\cdots , A_p] \cdot  C_Z(0) \cdot [ A_1 , A_2 , \hdots , A_p ]^\top \;.\label{eq:SigmafromYuleWalker}
\end{eqnarray}

\subsection{The time-reversed VAR($p$) process}
\label{app:timereversedVARp}

The results of Bartlett on the analytical description of time-reversed VAR(1) processes have been generalized to VAR($p$) processes by Andel in 1972 \cite{Andel1972}. Given a $2$-dimensional VAR($p$) process $\bz_t$ as in \eqref{eq:VARp}, Andel considers a second VAR($p$) process 
\begin{equation}
\boldsymbol{\mathfrak{z}}_t = \tilde{A}_1 \boldsymbol{\mathfrak{z}}_{t-1} + \tilde{A}_2 \boldsymbol{\mathfrak{z}}_{t-2} + \ldots + \tilde{A}_p \boldsymbol{\mathfrak{z}}_{t-p} + \mathfrak{e}_t \label{eq:VARprevAndel} \;,
\end{equation}
where $\mathfrak{e}_t$ is white noise with covariance matrix $\tilde{\Sigma} = \langle \mathfrak{e}_t \mathfrak{e}_t^\top \rangle$. 

Now, denote with $Q := C_Z(0)^{-1} $ the inverse of the covariance of $Z_t$, with block matrix notation 
\begin{eqnarray*}
\arraycolsep=1pt  
Q =: (Q_{lk})_{l,k=1}^p =  \begin{bmatrix} Q_{1,1} & Q_{1,2}  &\cdots & Q_{1,p} \\ Q_{2,1} & Q_{2,2}& \cdots & Q_{2,p}\\  \vdots &  \vdots & \ddots &  \vdots  \\ Q_{p,1}& Q_{p,p-1}& \cdots & Q_{p,p}\end{bmatrix} \in \R^{2p \times 2p} \;,
\end{eqnarray*} 
 where $Q_{lk}$ are $2 \times 2$ blocks.

Andel proves that $C_{\boldsymbol{\mathfrak{z}}}(h) = C_{\bz}(-h)$ for all $h\in \Z$  (that is, $\boldsymbol{\mathfrak{z}}_t$ has the same cross-covariance matrices as $\bz_t$ reversed in time), if and only if $\tilde{A}_1, \ldots , \tilde{A}_p$ and $\tilde{\Sigma}$ are defined as follows:  for $1\leq j \leq p$,
\begin{align}
\tilde{A}_j &= -(Q_{pp} + A_p^\top \Sigma^{-1} A_p)^{-1} (Q_{p,p-j} + A_p^\top \Sigma^{-1} A_{p-j}) \label{eq:AndelA} 
\intertext{and}
\tilde{\Sigma} &= (Q_{pp} + A_p^\top\Sigma^{-1} A_p)^{-1} \;,\label{eq:AndelSigma} 
\end{align}
where $Q_{p,0} := 0$ and $A_0 := -I$. Andel further proves that $\tilde{A}_p \neq 0$, if and only if $A_p \neq 0$
, and that, if $\bz_t$ is stable, so is $\boldsymbol{\mathfrak{z}}_t$. 

Note that, while we only treat bivariate VAR processes in this paper, the analytic description reviewed above holds for arbitrary dimensionality. 

\subsection{Proof that $\det(\Sigma) = \det (\tilde{\Sigma})$ for general p -- this completes the proof of Theorem \ref{prop:toprove}}
\label{app:proofequaldets}

Given Andel's result, we can complete the proof for Theorem \ref{prop:toprove}. The only missing part of the proof (cf. Section~\ref{sec:proofoftheoremVAR1}) is the proof of \eqref{eq:equaldets}, \mbox{$\det(\Sigma) = \det (\tilde{\Sigma})$}, for arbritrary $p \in \N \setminus \{ 0 \}$.

\emph{Preliminaries.} We will make use of the following well-known equalities.
Let $K$ be a positive definite matrix with \mbox{$L = K^{-1}$}, and let
\begin{equation*}
 K = \begin{bmatrix}
       K_{1,1} & K_{1,2} \\ K_{2,1} & K_{2,2}
      \end{bmatrix},
      \qquad
 L = \begin{bmatrix}
       L_{1,1} & L_{1,2} \\ L_{2,1} & L_{2,2}
      \end{bmatrix}
\end{equation*}
be the block matrix notations of $L$ and $K$, where $K_{1,1}$ is a square matrix of the same size as $L_{1,1}$.
Then the standard Schur complement formula (e.g.~\cite{Ouellette1981}, Theorem 2.7) is given as 
\begin{align}
 L_{2,2} = \left[K_{2,2} - K_{2,1}K_{1,1}^{-1}K_{1,2} \right]^{-1}. \label{eq:blockmatrixInverse}
\end{align}
Let $T$ and $W$ be invertible matrices, then for all matrices $U$ and $V$ of fitting size
\begin{equation}
 \det(T + UWV) = \det(W^{-1} + VT^{-1}U)\det(T)\det(W) \;.\label{eq:matrixDeterminantLemma}
\end{equation}
This relation is known as the generalized \emph{matrix determinant lemma} and a straightforward extension of Sylvester's determinant theorem~\eqref{eq:Sylvester}. 

Let $K$ be a matrix with block notation as above, and $K_{1,1}$ be invertible, then, (see e.g.~\cite{Ouellette1981}, Theorem 2.1)
\begin{equation}
 \det(K) = \det(K_{1,1})\det(K_{2,2} - K_{2,1}K_{1,1}^{-1}K_{1,2}) \;.\label{eq:blockmatrixDeterminantFactorization}
\end{equation}

Let us also introduce the following notation for the blocks of $C_Z(0)$:
 \setlength{\arraycolsep}{2pt} 
\begin{equation*}
C_Z(0) = \begin{bmatrix}
            C_{Z^{\setminus p}} & R^\top \\ R & C_\bz(0)
           \end{bmatrix} =\begin{bmatrix} C_\bz(0) & \bar{R} \\ \bar{R}^\top & C_{Z^{\setminus p}} \end{bmatrix} \;,
\end{equation*}
where we define
\begin{align*} 
R := & \begin{bmatrix} C_\bz(p-1)^\top & C_\bz (p-2)^\top & \cdots & C_\bz (1)^\top \end{bmatrix} \in \R^{2 \times 2(p-1)} \\
\bar{R} := &\begin{bmatrix}
               C_\bz(1) & \ldots & C_\bz(p-1)
              \end{bmatrix} \in \R^{2 \times 2(p-1)} \;,
\end{align*}
and
\begin{equation*} C_{Z^{\setminus p}}:= \begin{bmatrix} C_\bz(0)  &\cdots & C_\bz(p-2) \\   \vdots  & \ddots &  \vdots  \\ C_\bz(2-p)& \cdots & C_\bz(0) \end{bmatrix} \in \R^{2(p-1) \times 2(p-1)} \;.
\end{equation*}

\noindent \emph{Step 1: Analytic expression for $C_\bz(0) - \bar{R} C_{Z^{\setminus p}}^{-1}\bar{R}^\top$.}

We first prove that 
\begin{equation} 
C_\bz(0) - \bar{R} C_{Z^{\setminus p}}^{-1}\bar{R}^\top = \Sigma + A_pQ_{pp}^{-1} A_p^\top \;,
\label{eq:asterisk}
\end{equation}
that is, the residual variance when regressing $z_t$ on $z_{t-1}, \ldots, z_{t-(p-1)}$ given by $C_\bz(0) - \bar{R} C_{Z^{\setminus p}}^{-1}\bar{R}^\top$ can be expressed as the sum of $A_pQ_{pp}^{-1} A_p$ and the residual variance when regressing $z_t$ on $z_{t-1}, \ldots, z_{t-(p-1)},   z_{t-p}$, given by $\Sigma$.

Recall the Yule-Walker equation \eqref{eq:YuleWalkerVARp}
\begin{equation*}
C_Z(0) = \bold A \cdot C_Z(0) \cdot \bold A^\top + \Sigma_E \;.
\end{equation*}
and let us rewrite
\[\Sigma_E = \begin{bmatrix}
              \Sigma & 0 \\ 0 & 0
             \end{bmatrix}\quad \text{and}\quad
\bold A = \begin{bmatrix}
              \bold A_{\setminus p} & A_p \\ I & 0
             \end{bmatrix},
\]
where we define
\[\bold A_{\setminus p} := \begin{bmatrix} A_1 & \ldots & A_{p-1} \end{bmatrix} \in \R^{2 \times 2(p-1)} \;.\]
The  Yule-Walker equation can then be written in blocks as
\begin{equation*} 
\arraycolsep=2pt
\begin{bmatrix} C_\bz(0) & \bar{R} \\ \bar{R}^\top & C_{Z^{\setminus p}} \end{bmatrix} = \begin{bmatrix} \bold A_{\setminus p} & A_p \\ I & 0  \end{bmatrix} 
 		\begin{bmatrix} C_{Z^{\setminus p}} & R^\top \\ R & C_\bz(0) \end{bmatrix} 
		\begin{bmatrix}  \bold A_{\setminus p}^\top & I \\ A_p^\top & 0 \end{bmatrix} 
 	+ \begin{bmatrix}
              \Sigma & 0 \\ 0 & 0
             \end{bmatrix}  \;. 
\end{equation*}
We see from the top line that
\begin{align}
\bar{R} &= \bold A_{\setminus p}C_{Z^{\setminus p}} + A_pR\nonumber\\
\Leftrightarrow  \quad A_pR &= \bar{R} - \bold A_{\setminus p}C_{Z^{\setminus p}} \;,\label{eq:Rtilde}
\end{align}
and that
\setlength{\medmuskip}{2mu}
\begin{align}
C&_\bz(0) \nonumber \\
& =  \bold A_{\setminus p} C_{Z^{\setminus p}} \bold A_{\setminus p}^\top + \bold A_{\setminus p}R^\top A_p^\top + A_p  R \bold A_{\setminus p}^\top + A_pC_\bz(0)A_p^\top + \Sigma \nonumber\\
&\overset{\eqref{eq:Rtilde}}{=} - \bold A_{\setminus p} C_{Z^{\setminus p}} \bold A_{\setminus p}^\top + \bold A_{\setminus p}\bar{R}^\top + \bar{R}\bold A_{\setminus p}^\top + A_pC_\bz(0)A_p^\top + \Sigma \label{eq:Cz0Representation}
\end{align}
from which we conclude that
\begin{align*}
 \Sigma + &A_pQ_{pp}^{-1} A_p^\top \overset{\eqref{eq:blockmatrixInverse}}{=} \Sigma + A_p\left[C_\bz(0) - R C_{Z^{\setminus p}}^{-1}R^\top\right] A_p^\top\\
 &= \Sigma + A_pC_\bz(0)A_p^\top - A_pR C_{Z^{\setminus p}}^{-1}R^\top A_p^\top \\
 &\overset{\eqref{eq:Rtilde}}{=} \Sigma + A_pC_\bz(0)A_p^\top - \left[\bar{R} - \bold A_{\setminus p}C_{Z^{\setminus p}}\right]C_{Z^{\setminus p}}^{-1}\left[\bar{R} - \bold A_{\setminus p}C_{Z^{\setminus p}}\right]^\top\\
 &= \Sigma + A_pC_\bz(0)A_p^\top - \bar{R}C_{Z^{\setminus p}}^{-1}\bar{R}^\top + \bold A_{\setminus p}\bar{R}^\top + \bar{R}\bold A_{\setminus p}^\top - \bold A_{\setminus p}C_{Z^{\setminus p}}\bold A_{\setminus p}^\top\\
 &\overset{\eqref{eq:Cz0Representation}}{=} C_\bz(0) - \bar{R}C_{Z^{\setminus p}}^{-1}\bar{R}^\top \;.
\end{align*}

\pagebreak
\noindent \emph{Step 2. Derive $\det \Sigma = \det \tilde{\Sigma}$ from Andel  and \eqref{eq:asterisk}.}

From Andel~\eqref{eq:AndelSigma} we know that \mbox{$\tilde{\Sigma} = (Q_{pp} + A_p^\top\Sigma^{-1} A_p)^{-1}$}.
As $Q$ is positive definite, $Q_{pp}$ is invertible such that
\[\frac{1}{\det(\tilde{\Sigma})} = \det(Q_{pp} + A_p^\top\Sigma^{-1} A_p) \overset{\eqref{eq:matrixDeterminantLemma}}{=} \det(\Sigma + A_pQ_{pp}^{-1} A_p^\top)\frac{\det(Q_{pp})}{\det(\Sigma)} \;.\] It therefore suffices to show that
\begin{align}
\det(\Sigma + A_pQ_{pp}^{-1} A_p^\top) = \det(Q_{pp}^{-1}) \;.
\label{eq:lastthing}
\end{align}
Drawing on Step 1, Equation \eqref{eq:lastthing} can be proven as follows: 
 \setlength{\arraycolsep}{0pt} 
\begin{alignat*}{2}
 && \det(\Sigma + A_pQ_{pp}^{-1} A_p^\top) &= \det(Q_{pp}^{-1})\\
 & \overset{\eqref{eq:asterisk}, \eqref{eq:blockmatrixInverse}}{\Leftrightarrow} & \det(C_\bz(0) - \bar{R} C_{Z^{\setminus p}}^{-1}\bar{R}^\top) &= \det(C_\bz(0) - R C_{Z^{\setminus p}}^{-1}R^\top)\\
& \overset{\eqref{eq:blockmatrixDeterminantFactorization}}{\Leftrightarrow} &\det\left(\begin{bmatrix} C_{Z^{\setminus p}} & \bar{R}^\top \\ \bar{R} & C_\bz(0) \end{bmatrix}\right) \det(C_{Z^{\setminus p}}^{-1} )&= \det(C_Z(0)) \det(C_{Z^{\setminus p}}^{-1} )\\
& \Leftrightarrow &\det\left(\begin{bmatrix} C_{Z^{\setminus p}} & \bar{R}^\top \\ \bar{R} & C_\bz(0) \end{bmatrix}\right) &= \det(C_Z(0)) \;.
\end{alignat*}
Switching rows or columns of a matrix leaves its determinant invariant up to a factor of $(-1)^{i+j}$, where $i$ and $j$ are corresponding row or column indices. In the following, we perform block-wise rotation of a matrix block to the bottom, and to right, respectively. This is a concatenation of several row and column switches giving us a factor
$(-1)^r$ for a given $r$.  Note that $r$ is the same for both operations due to their symmetric behavior.
Therefore we have
\begin{align*}
 \det\left(\begin{bmatrix} C_{Z^{\setminus p}} & \bar{R}^\top \\ \bar{R} & C_\bz(0) \end{bmatrix}\right) &= (-1)^r \det\left(\begin{bmatrix} \bar{R} & C_\bz(0) \\ C_{Z^{\setminus p}} & \bar{R}^\top \end{bmatrix}\right)\\
 & = (-1)^{2r} \det\left(\begin{bmatrix} C_\bz(0) & \bar{R} \\ \bar{R}^\top & C_{Z^{\setminus p}}  \end{bmatrix}\right) \\
&= \det\left(\begin{bmatrix} C_\bz(0) & \bar{R} \\ \bar{R}^\top & C_{Z^{\setminus p}} \end{bmatrix}\right) = \det(C_Z(0)) \;,
\end{align*}
which completes the proof.

\subsection{Proof of Theorem 2}
\label{app:thm1InequalityDefiniteness}

\emph{Proof "$\Leftarrow$"} 

 If $A_1, \ldots , A_p$ are diagonal, then under assumption (A2) of uncorrelated residuals, $x_t$ and $y_t$ are independent. It follows immediately from the argumentation in Section~\ref{sec:NoiseRobustness} that $\tilde{D}_{x \rightarrow y}^{(net)} = 0$. 
\hfill \qed

\emph{Proof "$\Rightarrow$"}

Let  $\tilde{D}_{x \rightarrow y}^{(net)} = 0$.  From \eqref{eq:abschaetzung1} and \eqref{eq:abschaetzung2} in  Theorem~\ref{prop:toprove} we necessarily have $\tilde{\Sigma}_{xx} = \Sigma_{xx}$
and $\tilde{\Sigma}_{yy} = \Sigma_{yy}$.

From the equality of determinants  $\det(\Sigma) = \det(\tilde{\Sigma})$ in \eqref{eq:equaldets} and assumption (A2) of uncorrelated residuals, $\Sigma_{xy} = 0$, it follows:
\begin{equation*}
\Sigma = \tilde{\Sigma} = diag \, .
\label{eq:equalSigmas}
\end{equation*}

Now, the induction proof of the following statement completes the whole proof.

\noindent \emph{Proof via induction:}

\noindent \emph{Statement S(p):  In any stable bivariate VAR(p) process~\eqref{eq:VARp} with the time-reversed representation \eqref{eq:VARprev} fulfilling (A1), (A2) and (A3), and 
\begin{equation*}
(AS): \quad \Sigma = \tilde{\Sigma} = diag \;,
\end{equation*} 
the $A_1, \ldots , A_p$ are diagonal. }

\noindent \emph{Preliminaries}. 

First,  note that (AS)  implies that  the time-reversed coefficient matrices $\tilde{A}_1, \tilde{A}_2, \ldots, \tilde{A}_p$ are lower triangular, since
\begin{equation}
\tilde{\Sigma}_{xx} \overset{(AS)}{=}\Sigma_{xx} \overset{(A1)}{=} \Sigma_x \overset{\eqref{eq:uniequi}}{=} \tilde{\Sigma}_x \,.
\label{eq:Arevlower}
\end{equation}
From Andel \eqref{eq:AndelA} we know that 
\begin{equation*} 
\underbrace{\vphantom{A_p^\top}\tilde{A}_p}_{lower} \underbrace{\vphantom{A_p^\top}\Sigma}_{diag} = \underbrace{\vphantom{A_p^\top}\tilde{\Sigma}}_{diag}\underbrace{A_p^\top}_{upper} \, .
\end{equation*}
With (AS),
\begin{equation}
A_p = \tilde{A_p} = diag
\label{eq:equalAp}
\end{equation} 
immediately follows.

\noindent \emph{Basis: Show that the statement holds for $p = 1$}. 

Follows directly from \eqref{eq:equalAp}.

\noindent \emph{Inductive step: Show that if S(p-1) then S(p)}. 

From Andel \eqref{eq:AndelSigma}:
\begin{align}
 & \Sigma \overset{(AS)}{=} \tilde{\Sigma} \overset{\eqref{eq:AndelSigma} }{=} (Q_{pp} + A_p^\top\Sigma^{-1} A_p)^{-1} = (Q_{pp} + A_p^2\Sigma^{-1})^{-1} \nonumber \\
 \Rightarrow \; & \Sigma^{-1} = Q_{pp} + A_p^2\Sigma^{-1} \nonumber \\
 \Rightarrow \; & Q_{pp} = (I - A_p^2)\Sigma^{-1} = diag. \label{eq:Qppdiag}
\end{align}
Denote with $\tilde{Z}_t$ the VAR(1) representation of the time-reversed process $\tilde{\bz}_t$, and with $\tilde{Q} := C_{\tilde{Z}}(0)^{-1}$ its inverse covariance matrix with block notation $\tilde{Q} =: (\tilde{Q}_{lk})_{l,k=1}^p$.
Then, due to symmetry, 
\begin{equation}
\tilde{Q}_{pp} = (I - \tilde{A}_p^2)\tilde{\Sigma}^{-1} \overset{\eqref{eq:equalAp},(AS)}{=} (I - A_p^2)\Sigma^{-1} = Q_{pp}.
\label{eq:QppEqual}
\end{equation}
Let us now define the VAR($p-1$) process $\bz_t'$ by
\begin{equation}
\bz_t' = \sum_{i=1}^{p-1}B_i\bz_{t-i}' + \xi_t, \qquad \langle \xi_t \xi_t^\top \rangle = \Sigma',
\end{equation}
where 
\begin{align*}
[B_1, \ldots , B_{p-1}] &:= [C_\bz(1), \ldots , C_\bz(p-1)] \cdot C_{Z^{\setminus p}}(0)^{-1} \\
& = \bar{R} \cdot C_{Z^{\setminus p}}(0)^{-1} \\
\Sigma' &:= C_\bz(0) - \bar{R} C_{Z^{\setminus p}}(0)^{-1} \bar{R}^\top 
\end{align*}
arise from solving the Yule-Walker equations with respect to $C_\bz(0),C_\bz(1),\ldots,C_\bz(p-1)$. They are the least squares solution when regressing $\bz_{t}$ onto $\bz_{t-1},\ldots,\bz_{t-p+1}$.

Denote with $\tilde{B}_1, \ldots , \tilde{B}_{p-1}$ and $\tilde{\Sigma}'$ the time-reversed coefficients and residual covariance matrix of $\tilde{\bz}_t' := \bz_{-t}'$.
By definition of time reversal, the coefficients of $\tilde{\bz}_t'$ are the solution of the Yule-Walker equations with respect to $C_\bz(0),C_\bz(-1),\ldots,C_\bz(1-p)$. They are the least squares solution when regressing $\tilde{\bz}_{t}$ onto $\tilde{\bz}_{t-1},\ldots,\tilde{\bz}_{t-p+1}$. 

We now show that all assumptions in $S(p-1)$ hold for $\bz'_t$. 

\begin{itemize}[leftmargin=1cm]
\item[{[A0]}] \emph{Stability of $\bz_t'$.} \\ 
Since $\bz_t$ is stable, performing an insufficient lag-order fit (which is the case for $\bz_t'$) preserves stability (see \cite{whittle1963fitting}).
\item[{[A1]}] \emph{Lower triangularity of $B_1, \ldots, B_{p-1}$.} \\
Equations (8) and (10) in \cite{sowell1989decomposition} state the following relations: {\small
\begin{align} 
 \hspace{-19pt}\left[A_1 , \ldots , A_{p-1} \vphantom{\tilde{B}_{1}} \right] &= \left[B_1,\ldots,B_{p-1}\vphantom{\tilde{B}_{1}}\right] - A_p\left[\tilde{B}_{p-1},\ldots,\tilde{B}_{1}\right]\label{eq:Afactorization}\\
 \nonumber\hspace{-19pt}\left[ \tilde{A}_1 , \ldots , \tilde{A}_{p-1} \right] &= \left[\tilde{B}_1,\ldots,\tilde{B}_{p-1}\right] - \underbrace{\tilde{A}_p}_{A_p}\left[B_{p-1},\ldots,B_{1}\vphantom{\tilde{B}_{1}}\right] 
\end{align}}
Since all $A_i$ and $\tilde{A}_i$ are lower triangular ((A1), \eqref{eq:Arevlower}), and $A_p$ is diagonal, we deduce:
\[\begin{array}{ll} \left[\tilde{B}_{p-k}\right]_{12} = \left[A_p\right]_{11}^{-1}\left[B_{k}\right]_{12}\vspace{0.1cm} \\
\left[\tilde{B}_{p-k}\right]_{12} = \left[A_p\right]_{11}\left[B_{k}\right]_{12} \end{array}, \quad \forall k\in\{1,\ldots,p-1\} \;. \] 
This may only be fulfilled when either $\left[A_p\right]_{11} = 1$ or all $B_{k}$ (and $\tilde{B}_{k}$) are lower triangular.
Recalling that $A_p$ is diagonal and $Q_{pp} = (I - A_p^2)\Sigma^{-1}$ is invertible, necessarily we have $\left[A_p\right]_{11} \neq 1$. Therefore, $B_1, \ldots, B_{p-1}$ are lower triangular. 
\item[{[A2]}] \emph{Diagonality of $\Sigma'$} 
\item[{[AS]}]  \emph{and $\Sigma' = \tilde{\Sigma}'$.} \\
From \cite{sowell1989decomposition}, page 6, we have that 
\begin{equation*}
Q_{pp}^{-1} = \tilde{\Sigma}'   \quad \mbox{ and }  \tilde{Q}_{pp}^{-1} = \Sigma'.
\end{equation*}
From \eqref{eq:Qppdiag} and \eqref{eq:QppEqual}, 
\begin{equation*}
\Sigma' = \tilde{\Sigma}' = diag
\end{equation*}
immediately follows. 
\item[{[A3]}] \emph{Invertibility of $C_{\bz'}(0)$.} \\
By construction $\bz_t'$ is the solution of the Yule-Walker equations with respect to the autocovariances \mbox{$C_\bz(0),\ldots, C_\bz(p-1)$}. Thus, $C_{\bz'}(0) = C_\bz(0)$, which is invertible. 
\end{itemize}

\noindent Therefore, we can apply the induction claim $S(p-1)$ for $\bz_t'$ saying that $B_1, \ldots , B_{p-1}$ are diagonal. With the same argument $\tilde{B}_1, \ldots , \tilde{B}_{p-1}$ are also diagonal.

Using the factorization \eqref{eq:Afactorization}
\[A_k = \underbrace{B_k\vphantom{\tilde{B}_{1}}}_{diag} - \underbrace{A_p\vphantom{\tilde{B}_{1}}}_{diag}\underbrace{\tilde{B}_{p-k}}_{diag}, \quad k\in\{1,\ldots,p-1\}\]
we have that $A_1, \ldots, A_p$ are diagonal, which completes the proof. \qed

%
%

\ifCLASSOPTIONcaptionsoff
  \newpage
\fi

\bibliographystyle{IEEEtran}
\bibliography{causality_literature}

\end{document}